\documentclass[10pt]{article}
\usepackage[margin=2.5cm]{geometry}
\usepackage[utf8]{inputenc}
\usepackage{bussproofs}
\usepackage{mathrsfs}
\usepackage{amsmath}
\usepackage{amsthm}
\usepackage{amssymb}
\usepackage[dvipsnames]{xcolor}
\usepackage{stmaryrd}
\usepackage{calc}
\usepackage{cmll}
\usepackage{url}
\usepackage{amsmath}
\usepackage{longtable}
\usepackage{stmaryrd}
\usepackage{{xargs,ifthen,xstring}}
\usepackage{authblk}
\usepackage{txfonts}
\usepackage{array}
\newtheorem{theorem}{Theorem}
\newtheorem{convention}{Convention}
\newtheorem{definition}{Definition}
\newtheorem{proposition}{Proposition}

\newtheorem{corollary}{Corollary}

\title{A comparison of three kinds of monotonic proof-theoretic semantics and the base-incompleteness of intuitionistic logic}

\author{Antonio Piccolomini d'Aragona\\

Eberhard Karls Universit\"{a}t T\"{u}bingen, T\"{u}bingen, Germany\\

\texttt{antonio.piccolomini-daragona@uni-tuebingen.de}}
\date{}

\begin{document}

\maketitle

\begin{abstract}
    I deal with two approaches to proof-theoretic semantics: one based on argument structures and justifications, which I call \emph{reducibility semantics}, and one based on consequence among (sets of) formulas over atomic bases, called \emph{base semantics}. The latter splits in turn into a \emph{standard reading}, and a variant of it put forward by Sandqvist. I prove some results which, when suitable conditions are met, permit one to shift from one approach to the other, and I draw some of the consequences of these results relative to the issue of completeness of (recursive) logical systems with respect to proof-theoretic notions of validity. This will lead me to focus on a notion of base-completeness, which I will discuss with reference to known completeness results for intuitionistic logic. The general interest of the proposed approach stems from the fact that reducibility semantics can be understood as a labelling of base semantics with proof-objects typed on (sets of) formulas for which a base semantics consequence relation holds, and which witness this very fact. Vice versa, base semantics can be understood as a type-abstraction of a reducibility semantics consequence relation obtained by removing the witness of the fact that this relation holds, and by just focusing on the input and output type of the relevant proof-object.
\end{abstract}


\section{Introduction}

The name \emph{proof-theoretic semantics} (PTS for short) indicates a family of constructive semantics whose core notion for explaining meaning and for defining the notions of (logical) validity is not that of truth, as happens in model theory, but that of proof; the meaning of the non-logical vocabulary is accordingly fixed, not via model-theoretic mappings from the language onto suitable (typically set-theoretic) structures, but through sets of (sets of) rules governing deduction at the atomic level---see \cite{schroeder-heisterSE} for an overview.

PTS stems from Prawitz's work in proof theory, in particular from Prawitz's normalisation theorems for Gentzen's Natural Deduction \cite{gentzen, prawitz1965}. The first version of PTS is due to Prawitz himself; in this formulation, PTS is based on \emph{argument structures} and \emph{reductions} \cite{prawitz1971, prawitz1973}, so I will call it \emph{reducibility semantics}. Later on---maybe starting from the influential \cite{schroeder-heister2006}---argument structures were left aside, and PTS became a theory of consequence for formulas over sets of (sets of) atomic rules; the constructivist burden was put entirely on such sets, and this is why---partly following \cite{sandqvist}---I will indicate this approach as \emph{standard base semantics}. A variant of base semantics was provided by Sandqvist \cite{sandqvist}, whence I shall call it \emph{Sanqvist's base semantics}. It differs from the standard reading---and from Prawitz's original approach---in that it deals with disjunction in an elimination-like way, rather than in an introduction-based fashion.

Both reducibility semantics and base semantics (in its two variants) are expected to be semantics for constructive logics, and many completeness and incompleteness results have been obtained so far---see \cite{piecha} for an overview, whereas more recent results can be found in \cite{piechaschroeder-heister2019, stafford1, stafford2}. Here, I prove some general results on the relation between reducibility semantics and base semantics, and use them to obtain further insights on completeness and incompleteness issues. In particular, I shall be concerned by the conditions or the extent under which one is allowed to go from a consequence relation (possibly relative to specific formulas and over a specific ‘‘model") in one of the proof-theoretic semantics versions at issue here, to the same consequence in another version.

The interest of comparing the three approaches is twofold. As regards the relation between reducibility semantics and base semantics in the standard reading, both are ultimately based on the idea that the meaning of a logical constant $\kappa$ should be given by the conditions for \emph{introducing} $\kappa$ in formulas---thus coping with BHK semantics \cite{troelstravandalen} and Gentzen's semantic insight about introduction rules of Natural Deduction \cite{gentzen}. However, as already said above, Prawitz's picture is at first sight richer, since there the constructivist spirit is expressed, not only through the kind of ‘‘models" which the notion of consequence is defined over, but additionally through ‘‘witnesses" for the consequence relation itself---i.e., $A$ being a consequence of $\Gamma$ is \emph{defined} as existence of a \emph{valid argument} from $\Gamma$ to $A$. As also remarked above, standard base semantics instead drops valid arguments out, and one might at that point wonder whether the original constructivist spirit of Prawitz's approach is respected. If the consequence relation is no longer ‘‘witnessed" by suitable proof-objects, we are left with ‘‘models" of a special kind only, so we may ask whether some ‘‘degrees of constructiveness" got lost in the pruning. I shall prove below that, to a \emph{large} extent, nothing is lost, and hence that Prawitz's semantics, when some constraints are met, is in fact equivalent to standard base semantics. It follows that all the (mostly in)completeness results which hold for the former apply to a certain understanding of the latter as well. This is admittedly not a major result, as it can be looked at as nothing but a ‘‘decoration" of the general theorems established in \cite{piechaschroeder-heister2019}, and was stated by \cite{piccolominibook} in the context of Prawitz's \emph{theory of grounds} \cite{prawitz2015}---for the non-monotonic approach, but easily extendable to the monotonic one. However, this ‘‘decoration" involves steps which show that, in the specific framework of reducibility semantics, a detailed proof involves some not-so-trivial aspects.

The second point of interest touches upon the relation between Prawitz's reducibility semantics and Sanqvist's base semantics---hence, given the equivalence between the former and base semantics in the standard reading, also between standard base semantics and its Sandqvist variant. Since Sandqvist's base semantics does without argument structures and reductions too, one may have also here issues about dropping out the ‘‘witnesses" of consequence that one had in Prawitz. But there is more than that. The elimination-based approach to disjunction employed by Sandqvist induces a structural difference with respect to Prawitz's introduction-based approach, which triggers intuitionistic completeness in some most relevant cases, \emph{contra} the provable intuitionistic incompleteness that we have for (base semantics in the standard reading, hence for) Prawitz's reducibility semantics (see below for details). Given this (and given, incidentally, the acknowledged \emph{harmony} between introduction and elimination rules in Natural Deduction), it seems thus to be worth exploring to what extent Prawitz's reducibility semantics and base semantics in Sandqvist's reading can be connected to each other. Below, I shall prove that the two approaches \emph{can} be compared, although only on a ‘‘global" scale, and that this, together with the completeness-incompleteness mismatch, implies that they \emph{are not} comparable at the level of ‘‘models"---meaning that the models verifying certain pairs $(\Gamma, A)$ in one approach are not always models verifying the same pairs in the other, and vice versa. The ‘‘global" comparability also provides a sufficient condition for equivalence to hold between Prawitz's and Sandqvist's pictures at the level of logical validity. This may be of interest since, e.g., while intuitionistic logic ($\texttt{IL}$ for short) is incomplete over standard base semantics (and hence, via the equivalence mentioned above, over reducibility semantics too) relative to ‘‘models" of \emph{any} kind (more precisely, as we shall see below, relative to atomic bases of either limited or unlimited complexity), $\texttt{IL}$ is instead complete over Sandqvist's base semantics relative to ‘‘models" of a specific kind.

The result about the sufficient condition for equivalence between Prawitz's and Sandqvist's approaches at the level of logical validity requires introducing a notion of ‘‘point-wise" soundness and completeness of given logics $\Sigma$ relative to proof-theoretic semantics. Here ‘‘point-wise" means, roughly, that validity over a model is implied (soundness) or implies (completeness) derivability in $\Sigma$ plus the model (a notion which make sense given, as we shall see, the ‘‘deductive" nature of models in proof-theoretic semantics). Besides illuminating the relation between Prawitz's and Sandqvist's frameworks, the notions of ‘‘point-wise" soundness and completeness will be shown to be of interest in themselves. In this connection, I shall prove a result of ‘‘point-wise" incompleteness (with respect to all the three proof-theoretic frameworks at issue here) for a class of super-intuitionistic logics, which includes $\texttt{IL}$.

The structure of the paper is as follows. By limiting myself to a propositional language, I start with an overview of atomic rules (Section 2). Then I provide an outline of base semantics and reducibility semantics (Section 3). In Section 4, I prove the general results mentioned above and draw some consequences from them, both relative to the general issue of completeness of given (recursive) systems, and relative to a notion of base-completeness (Section 5).

\section{Language and atomic bases}

\begin{definition}
    The \emph{language} $\mathscr{L}$ is given by the grammar

    \begin{center}
        $X \coloneqq p, q, r, s, t, ...$ (infinitely many) $ \ | \ \bot \ | \ X \wedge X \ | \ X \vee X \ | \ X \rightarrow X$
    \end{center}
    where $\bot$ is a constant atom for absurdity and $\neg X \coloneqq X \rightarrow \bot$.
\end{definition}
\noindent I will use the following notation:
\begin{itemize}
    \item $\texttt{ATOM}_\mathscr{L}$ is the set of the atomic formulas of $\mathscr{L}$, i.e., $\{p, q, r, s, t, ...\} \cup \{\bot\}$;
    \item $\texttt{FORM}_\mathscr{L}$ is the set of the formulas of $\mathscr{L}$;
    \item $A, B, C, ...$ indicate arbitrary formulas;
    \item $\Gamma, \Delta, \Theta, ...$ indicate arbitrary sets of formulas. It is important to remark that the latter will be \emph{always} assumed to be \emph{finite}. This is to avoid some issues, pointed out by \cite{stafford1}, concerning properties of compactness, monotonicity, and compact monotonicity of the consequence relations in the comparison of reducibility semantics and base semantics. The limitation to finite sets of formulas will be on the other hand sufficient for raising my points---in particular, for transferring completeness or incompleteness results from one approach to the other.
\end{itemize}

\noindent I now define the notion of atomic base over $\mathscr{L}$. Atomic bases, however, require a preliminary definition of the notion of atomic rule over $\mathscr{L}$. The definition of atomic rules is by induction on what \cite{piechaschroeder-heisterdecampossanz} calls the level of atomic rules.

\begin{definition}
    Any atom is an \emph{atomic rule of level $0$}. An \emph{atomic rule of level 1} is
    \begin{prooftree}
        \AxiomC{$A_1$}
        \AxiomC{$\dots$}
        \AxiomC{$A_n$}
        \TrinaryInfC{$A$}
    \end{prooftree}
    with $A_1, ..., A_n, A \in \texttt{ATOM}_\mathscr{L}$. Given sets of atomic rules $\Re_1, ..., \Re_n$ whose maximal level is $\kappa - 1 \geq 0$, where each $\Re_i$ may be empty ($i \leq n$), we say that

\begin{prooftree}
    \AxiomC{$[\Re_1]$}
    \noLine
    \UnaryInfC{$A_1$}
    \AxiomC{}
    \noLine
    \UnaryInfC{$\dots$}
    \AxiomC{$[\Re_n]$}
    \noLine
    \UnaryInfC{$A_n$}
    \TrinaryInfC{$A$}
\end{prooftree}
with $A_1, ..., A_n, A \in \texttt{ATOM}_\mathscr{L}$ is an \emph{atomic rule of level $\kappa + 1$}.
\end{definition}

\noindent Square brackets in Definition 2 indicate that the rule discharges lower level atomic rules in its premises. Thus, an atomic rule of level 3 has the form

\begin{prooftree}
    \AxiomC{$[R_{1, 1}]$}
    \AxiomC{}
    \noLine
    \UnaryInfC{$\dots$}
    \AxiomC{$[R_{1, m_1}]$}
    \noLine
    \TrinaryInfC{$A_1$}
    \AxiomC{}
    \noLine
    \UnaryInfC{$\dots$}
    \AxiomC{$[R_{n, 1}]$}
    \AxiomC{}
    \noLine
    \UnaryInfC{$\dots$}
    \AxiomC{$[R_{n, m_n}]$}
    \noLine
    \TrinaryInfC{$A_n$}
    \TrinaryInfC{$A$}
\end{prooftree}
where, for every $i \leq n, j \leq m_i$, either $R_{i, j} = B \in \texttt{ATOM}_\mathscr{L}$, or

\begin{prooftree}
    \AxiomC{$B_1$}
    \AxiomC{$\dots$}
    \AxiomC{$B_s$}
    \LeftLabel{$R_{i, j} = \ $}
    \TrinaryInfC{$B$}
\end{prooftree}
with $B_1, ..., B_s, B \in \texttt{ATOM}_\mathscr{L}$.

\noindent I will adopt the convention that atomic bases are---borrowing the terminology of \cite{piechaschroeder-heisterdecampossanz}---of an \emph{intuitionistic} kind, namely, that they contain a rule of \emph{atomic explosion} for every atom in the language.

\begin{convention}
    Every $\mathfrak{B}$ contains the rules
    \begin{prooftree}
        \AxiomC{$\bot$}
        \RightLabel{$\texttt{AtExp}$}
        \UnaryInfC{$A$}
    \end{prooftree}
    for every $A \in \texttt{ATOM}_\mathscr{L}$.
\end{convention}

\begin{definition}
    An \emph{atomic base of level $n$} is a set of atomic rules $\mathfrak{B} = \{\texttt{AtExp}, R_1, ..., R_m\}$ with $\texttt{max}\{\mathfrak{L}(R_1), ..., \mathfrak{L}(R_m)\} = n$, where $\mathfrak{L}(R_i)$ indicates the level of $R_i$ ($i \leq m$).
\end{definition}
\noindent Some remarks on atomic bases are now in order. First, the fact that I used a finite index $m$ to refer to the rules $R_1, ..., R_m$ of $\mathfrak{B}$, \emph{does not imply} that $\mathfrak{B}$ should be understood as \emph{finite}. It just means that---as done with $\texttt{AtExp}$---each $R_i \ (i \leq m)$ is thought of as given schematically---where schemes are allowed to have infinitely many instances. Thus, atomic bases are allowed to contain infinitely many atomic rules. An example, which I shall also use below, is $\{R\}$ where $R$ is

\begin{prooftree}
    \AxiomC{$A$}
    \AxiomC{$[B]$}
    \noLine
    \UnaryInfC{$D$}
    \AxiomC{$[C]$}
    \noLine
    \UnaryInfC{$D$}
    \TrinaryInfC{$D$}
\end{prooftree}
for some $A, B, C \in \texttt{ATOM}_\mathscr{L}$ and \emph{every} $D \in \texttt{ATOM}_\mathscr{L}$. Second, I shall use the following notation: $\mathbb{B}^n = \{\mathfrak{B} \ | \ \mathfrak{B}$ has level $m \leq n\}$. Third, the empty atomic base, written $\mathfrak{B}^\emptyset$, is the atomic base which only contains $\texttt{AtExp}$. Observe that, as per Definition 3, the latter is disregarded when counting the level of the atomic base, so $\mathfrak{B}^\emptyset$ has level $0$.

The requirement that atomic bases always contain $\texttt{AtExp}$ might be motivated as follows---see, e.g., \cite{piechaschroeder-heisterdecampossanz}. Suppose we explain the meaning of $\bot$ by saying that there is no $\mathfrak{B}$ such that $\bot$ is valid on $\mathfrak{B}$---in one of the senses of validity to be specified below. Then, for every $\mathfrak{B}$ and every $A \in \texttt{ATOM}_\mathscr{L}$, $\neg \neg A$ is valid on $\mathfrak{B}$---i.e., $\neg \neg A$ is logically valid. This is because, if we assume that $\neg A$ is valid on $\mathfrak{B}$, for every extension of $\mathfrak{B}$ where $A$ is valid, also $\bot$ is valid. On the other hand, we can always extend $\mathfrak{B}$ by adding $A$ as an axiom to it, whence we would have that $\bot$ is valid on this extension of $\mathfrak{B}$, which is ruled out by our semantic clause for $\bot$. The strategy I am adopting here, where $\bot$ is an atomic constant and atomic bases always contain $\texttt{AtExp}$, is one of the two major approaches for monotonic PTS, the other being the one where the meaning of $\bot$, now a nullary connective, is explained by requiring that $\bot$ is valid on $\mathfrak{B}$ if and only if every $A \in \texttt{ATOM}_\mathscr{L}$ is valid on $\mathfrak{B}$. In Proposition 2 below, I show that this distinction is harmless in the context of this paper. Recently, Barroso Nascimento, Pereira \& Pimentel have pointed out that another solution is possible, where atomic bases are forced to be consistent, thereby indirectly restricting the classes of extensions of the bases themselves---and ruling out the possibility of always adding an atom as an axiom to a given base \cite{barrosopereiranascimento}. I shall not take this alternative into account here but---also based on the results proved in \cite{barrosopereiranascimento}---it is my impression that the results proved in the present paper also apply to the variant in question.

\begin{definition}
    The \emph{derivations-set} $\texttt{DER}_\mathfrak{B}$ of $\mathfrak{B}$ is defined inductively as follows:

    \begin{itemize}
        \item any $A \in \texttt{ATOM}_\mathscr{L}$ is a single-node derivation in $\texttt{DER}_\mathfrak{B}$. It is a derivation of $A$ from $\emptyset$ if the node applies an axiom $A \in \mathfrak{B}$, or a derivation of $A$ from $A$ if $A$ is assumed as an atomic rule (whether or not $A \in \mathfrak{B}$);
        \item if the following are derivations in $\texttt{DER}_\mathfrak{B}$,
        \begin{prooftree}
            \AxiomC{$\mathfrak{C}_i, \Re_i$}
            \noLine
            \UnaryInfC{$\mathscr{D}_i$}
            \noLine
            \UnaryInfC{$A_i$}
        \end{prooftree}
        ($i \leq n$) where $\mathfrak{C}_i$ and $\Re_i$ are sets of atomic rules and $A_i \in \texttt{ATOM}_\mathscr{L}$ is the premise of an atomic rule $R$ of the form
        \begin{prooftree}
            \AxiomC{$[\Re_1]$}
            \noLine
            \UnaryInfC{$A_1$}
            \AxiomC{$\dots$}
            \AxiomC{$[\Re_n]$}
            \noLine
            \UnaryInfC{$A_n$}
            \RightLabel{$R$}
            \TrinaryInfC{$A$}
        \end{prooftree}
        then
        \begin{prooftree}
            \AxiomC{$\mathfrak{C}_1, [\Re_1]$}
            \noLine
            \UnaryInfC{$\mathscr{D}_1$}
            \noLine
            \UnaryInfC{$A_1$}
            \AxiomC{$\dots$}
            \AxiomC{$\mathfrak{C}_n, [\Re_n]$}
            \noLine
            \UnaryInfC{$\mathscr{D}_n$}
            \noLine
            \UnaryInfC{$A_n$}
            \RightLabel{$R$}
            \TrinaryInfC{$A$}
        \end{prooftree}
        is a derivation of $A$ from $\mathfrak{C}_1, ..., \mathfrak{C}_n$ in $\texttt{DER}_\mathfrak{B}$ if $R \in \mathfrak{B}$, or a derivation of $A$ from $\mathfrak{C}_1, ..., \mathfrak{C}_n, R$ if $R \notin \mathfrak{B}$.
    \end{itemize}
\end{definition}

\begin{definition}
    $A$ is \emph{derivable from $\Re$ in $\mathfrak{B}$}---written $\Re \vdash_\mathfrak{B} A$---iff there is $\mathscr{D} \in \texttt{DER}_\mathfrak{B}$ from $\Re$ to $A$.
\end{definition}

\begin{definition}
    $\mathfrak{C}$ is an \emph{extension} of $\mathfrak{B}$ iff $\mathfrak{B} \subseteq \mathfrak{C}$.
\end{definition}

\noindent Let me give an example of atomic derivation (in rules, left label indicates the name of the rule, with square brackets for indicating dischargement, while right label indicates that dischargement takes place at that rule). Suppose that our base $\mathfrak{B}$ contains one rule of level $0$, one rule of level $1$, one rule of level $2$ and one rule of level $3$, i.e.,

\begin{prooftree}
    \AxiomC{}
    \UnaryInfC{$p$}
    \AxiomC{$p$}
    \UnaryInfC{$q$}
    \AxiomC{}
    \LeftLabel{$[R_1]$}
    \UnaryInfC{$s$}
    \noLine
    \UnaryInfC{$t$}
    \AxiomC{$r$}
    \RightLabel{$R_1$}
    \BinaryInfC{$w$}
    \AxiomC{$p$}
    \AxiomC{$w$}
    \LeftLabel{$[R_2]$}
    \BinaryInfC{$t$}
    \noLine
    \UnaryInfC{$z$}
    \RightLabel{$R_2$}
    \UnaryInfC{$y$}
    \noLine
    \QuaternaryInfC{}
\end{prooftree}
and consider the atomic rules of level $0$ and $1$
\begin{prooftree}
    \AxiomC{}
    \LeftLabel{$R_3$}
    \UnaryInfC{$q$}
    \AxiomC{$t$}
    \AxiomC{$q$}
    \LeftLabel{$R_4$}
    \BinaryInfC{$z$}
    \AxiomC{$s$}
    \LeftLabel{$R_5$}
    \UnaryInfC{$t$}
    \AxiomC{$q$}
    \LeftLabel{$R_6$}
    \UnaryInfC{$r$}
    \noLine
    \QuaternaryInfC{}
\end{prooftree}
Then the following proves $R_3, R_4, R_5, R_6 \vdash_\mathfrak{B} y$
\begin{prooftree}
    \AxiomC{}
    \UnaryInfC{$p$}
    \AxiomC{}
    \LeftLabel{$[R_1]$}
    \UnaryInfC{$s$}
    \LeftLabel{$R_5$}
    \UnaryInfC{$t$}
    \AxiomC{}
    \UnaryInfC{$p$}
    \UnaryInfC{$q$}
    \LeftLabel{$R_6$}
    \UnaryInfC{$r$}
    \RightLabel{$R_1$}
    \BinaryInfC{$w$}
    \LeftLabel{$[R_2]$}
    \BinaryInfC{$t$}
    \AxiomC{}
    \LeftLabel{$R_3$}
    \UnaryInfC{$q$}
    \LeftLabel{$R_4$}
    \BinaryInfC{$z$}
    \RightLabel{$R_2$}
    \UnaryInfC{$y$}
\end{prooftree}

\noindent I use the following notation: $\mathfrak{C} \supseteq_n \mathfrak{B}$ indicates that $\mathfrak{C} \supseteq \mathfrak{B}$ and $\mathfrak{C} \in \mathbb{B}^n$.

\begin{proposition}
    $\{\mathfrak{B} \ | \ \mathfrak{B} \supseteq_n \mathfrak{B}^\emptyset\} = \mathbb{B}^n$.
\end{proposition}

\section{Base semantics and reducibility semantics}

As I am restricting myself to propositional logic, I will occasionally allow myself to indicate the universal and existential meta-quantifiers by the usual object-linguistic notations $\forall$ and $\exists$.

\subsection{Base semantics}

\begin{definition}
   That $A$ is a \emph{consequence of $\Gamma$ on $\mathfrak{B}$ of level $n$ in base semantics} is indicated by $\Gamma \models_{\mathfrak{B}, \ n} A$. It holds iff $\mathfrak{B} \in \mathbb{B}^n$ and
   \begin{itemize}
       \item[1.] $\Gamma = \emptyset \Longrightarrow$
       \begin{itemize}
           \item[(a)] $A \in \texttt{ATOM}_\mathscr{L} \Longrightarrow \ \vdash_\mathfrak{B} A$;
           \item[(b)] $A = B \wedge C \Longrightarrow \ \models_{\mathfrak{B}, \ n} B$ and $\models_{\mathfrak{B}, \ n} C$;
           \item[(c)] $A = B \vee C \Longrightarrow \ \models_{\mathfrak{B}, \ n} B$ or $\models_{\mathfrak{B}, \ n} C$;
           \item[(d)] $A = B \rightarrow C \Longrightarrow B \models_{\mathfrak{B}, \ n} C$;
       \end{itemize}
       \item[2.] $\Gamma \neq \emptyset \Longrightarrow \forall \mathfrak{C} \supseteq_n \mathfrak{B} \ (\models_{\mathfrak{C}, \ n} \Gamma \Longrightarrow \ \models_{\mathfrak{C}, \ n} A)$
   \end{itemize}
   where $\models_{\mathfrak{C}, \ n} \Gamma$ means that $\forall B \in \Gamma \ (\models_{\mathfrak{C}, \ n} B)$.
\end{definition}

\noindent A variant of base semantics was introduced by Sandqvist \cite{sandqvist}. In the present framework, it runs as follows.

\begin{definition}
    That $A$ is a \emph{consequence of} $\Gamma$ \emph{on} $\mathfrak{B}$ \emph{of level} $n$ \emph{in base semantics in Sandqvist sense} is indicated by $\Gamma \models^s_{\mathfrak{B}, \ n} A$. It holds iff $\mathfrak{B} \in \mathbb{B}^n$ and
    \begin{itemize}
        \item[1.] $\Gamma = \emptyset \Longrightarrow$
        \begin{itemize}
            \item[(a)] $A \in \texttt{ATOM}_\mathscr{L} \Longrightarrow \ \vdash_{\mathfrak{B}} A$;
            \item[(b)] $A = B \wedge C \Longrightarrow \ \models^s_{\mathfrak{B}, \ n} A$ and $\models^s_{\mathfrak{B}, \ n} B$;
            \item[(c)] $A = B \vee C \Longrightarrow \forall \mathfrak{C} \supseteq_n \mathfrak{B} \ \forall D \in \texttt{ATOM}_{\mathscr{L}} \ (B \models^s_{\mathfrak{C}, \ n} D$ and $C \models^s_{\mathfrak{C}, \ n} D \Longrightarrow \ \models^s_{\mathfrak{C}, \ n} D)$;
            \item[(d)] $A = B \rightarrow C \Longrightarrow B \models^s_{\mathfrak{B}, \ n} C$;
        \end{itemize}
        \item[2.] $\Gamma \neq \emptyset \Longrightarrow \ \forall \mathfrak{C} \supseteq_n \mathfrak{B} \ (\models^s_{\mathfrak{C}, \ n} \Gamma \Longrightarrow \ \models^s_{\mathfrak{C}, \ n} A)$
    \end{itemize}
\end{definition}

\noindent  We remark that, in Sandqvist, $\bot$ is understood as a \emph{non-atomic} constant, so we should add a new sign to our language, say $\bot^*$, as distinguished from $\bot$. However, over intuitionistic bases, $\bot^*$ and $\bot$ are equivalent \cite[but see also Proposition 2 below]{piechaschroeder-heisterdecampossanz}. In what follows, $\Vdash_n$ means $\models_n$ or $\models^s_n$---possibly with an index for atomic bases.

\begin{definition}
    That $A$ is a \emph{logical consequence of} $\Gamma$ \emph{of level} $n$ \emph{in base semantics (in Sandqvist sense)} is indicated by $\Gamma \Vdash_n A$. It holds iff $\forall \mathfrak{B} \in \mathbb{B}^n \ (\Gamma \Vdash_{\mathfrak{B}, \ n} A)$.
\end{definition}

\begin{proposition}
    $\Vdash_{\mathfrak{B}, \ n} \bot \Longleftrightarrow \ \forall A \in \texttt{ATOM}_\mathscr{L} \ (\Vdash_{\mathfrak{B}, \ n} A)$.
\end{proposition}

\begin{proof}
    ($\Longleftarrow$) is trivial since $\bot \in \texttt{ATOM}_\mathscr{L}$. ($\Longrightarrow$) $\Vdash_{\mathfrak{B}, \ n} \bot$ is equivalent to $\vdash_\mathfrak{B} \bot$ which, since $\mathfrak{B}$ contains $\texttt{AtExp}$, implies $\vdash_\mathfrak{B} A$ for every $A \in \texttt{ATOM}_\mathscr{L}$, which is again equivalent to $\Vdash_{\mathfrak{B}, \ n} A$.
\end{proof}

\begin{proposition}[Monotonicity of $\Vdash$]
    $\Gamma \Vdash_{\mathfrak{B}, \ n} A \Longleftrightarrow \forall \mathfrak{C} \supseteq_n \mathfrak{B} \ (\Gamma \Vdash_{\mathfrak{C}, \ n} A)$.
\end{proposition}

\begin{proof}
    ($\Longleftarrow$) is trivial by putting $\mathfrak{C} = \mathfrak{B}$. ($\Longrightarrow$) $\Gamma \Vdash_{\mathfrak{B}, \ n} A$ is equivalent to $\forall \mathfrak{C} \supseteq_n \mathfrak{B} \ (\Vdash_{\mathfrak{C}, \ n} \Gamma \Longrightarrow \ \Vdash_{\mathfrak{C}, \ n} A)$. The result follows by observing that any $\mathfrak{D} \supseteq_n \mathfrak{C}$ is also such that $\mathfrak{D} \supseteq_n \mathfrak{B}$.
\end{proof}

\begin{proposition}
    $\Gamma \Vdash_n A \Longleftrightarrow \Gamma \Vdash_{\mathfrak{B}^\emptyset, \ n} A$.
\end{proposition}

\begin{proof}
    It follows immediately from Proposition 1 and Proposition 3.
\end{proof}

\begin{corollary}
    $\Gamma \Vdash_n A \Longleftrightarrow \forall \mathfrak{B} \in \mathbb{B}^n \ (\Vdash_{\mathfrak{B}, \ n} \Gamma \Longrightarrow \ \ \Vdash_{\mathfrak{B}, \ n} A)$.
\end{corollary}

\begin{proof}
    By Proposition 4, $\Gamma \Vdash_n A$ is equivalent to $\Gamma \Vdash_{\mathfrak{B}^\emptyset, \ n} A$ which in turn, by points 2 in Definition 7 or 8, is equivalent to $\forall \mathfrak{B} \supseteq_n \mathfrak{B}^\emptyset \ (\Vdash_{\mathfrak{B}, \ n} \Gamma \Longrightarrow \ \Vdash_{\mathfrak{B}, \ n} A)$. By Proposition 1, this is equivalent to $\forall \mathfrak{B} \in \mathbb{B}^n \ (\Vdash_{\mathfrak{B}, \ n} \Gamma \Longrightarrow \ \Vdash_{\mathfrak{B}, \ n} A)$.
\end{proof}

\subsection{Reducibility semantics}

As said, reducibility semantics is instead based on the notions of argument structure and reduction. The latter are then used for defining the notions of argumental validity over a base, and of argumental validity in general. The presentation below is mainly inspired by \cite{prawitz1965, prawitz1973}, and partly by \cite{schroederheisternaturalext}---there are however some differences with respect to each of these sources, e.g., \cite{prawitz1973} does not use higher-level atomic rules, but just atomic rules of level at most $1$, and does not consider extensions of atomic bases. The kind of reducibility semantics that I put forward here, however, is to my mind better suited for the joint comparison with the two versions base-extension semantics which I discussed in the previous section---while hopefully in line with Prawitz's original picture.

\begin{definition}
    An \emph{argument structure over $\mathscr{L}$} is a pair $\langle T, \langle f, h, g \rangle \rangle$ such that
    
    \begin{itemize} 
    \item $T$ is a finite rooted tree (let $\prec$ denote the order relation over $T$) whose nodes are labelled by formulas of $\mathscr{L}$. Let us assume that the labels of the top-nodes of $T$ are partitioned into two groups:
    \begin{itemize}
    \item assumption-labels, which I indicate by $\texttt{as}(T)$, and
    \item axiom-labels, which I indicate by $\texttt{ax}(T)$;
    \end{itemize}
    \item $f$ is a function defined on some $\Gamma \subseteq \texttt{as}(T)$ such that, $\forall \mu \in \Gamma$, $\mu \prec f(\mu)$ and each $\mu \in \Gamma$ is labelled by an atom;
    \item $h$ is a function defined on some $H \subseteq \texttt{ax}(T)$ and is such that, $\forall \mu \in H$, it holds that $\mu \prec h(\mu)$, $h(\mu)$ and all its children are labelled by atoms, and there is no $\nu \in \Gamma$ such that $f(\nu) = h(\mu)$;
    \item $g$ is a function defined on some $J \subseteq \mathcal{P}(\prec)$ such that, for every $K \in J$,
    \begin{itemize}
    \item $K$ contains all and only the edges that link a given node $\mu_K$ to all its children, and
    \item both $\mu_K$ and its children are labelled by atoms, and
    \item there is no $\nu \in \Gamma$ such that $f(\nu) = \mu_K$, and
    \item the function is such that $\forall K \in J$, it holds that $\mu_K \prec g(K)$, $g(K)$ and all its children are labelled by atoms, and there is no $\xi \in \Gamma$ such that $f(\xi) = g(K)$.
    \end{itemize}
    \end{itemize}
\end{definition}

\begin{definition}
    Given $\mathscr{D} = \langle T, \langle f, h, g \rangle \rangle$ with $\texttt{as}(T) = \Gamma$ and root $A$, the elements of $\Gamma$ are the \emph{assumptions} of $\mathscr{D}$ and $A$ is the \emph{conclusion} of $\mathscr{D}$.
\end{definition}

\noindent The intended meaning of $\texttt{ax}(T)$ in Definition 10 is that the atoms which label the given top-nodes in $T$ are not assumptions, but level-$0$ rules. The meaning of $f, h, g$ in the same definition is as follows. In Natural Deduction terminology, they are discharge functions---see \cite{prawitz1965, schroederheisternaturalext, schroeder-heister2006}. Top-nodes $\mu \in \texttt{as}(T)$ in the domain of $f$ are assumptions discharged throughout $T$, while top-nodes $\mu \in \texttt{ax}(T)$ and sets of edges $\mu \in J$ in the domain of $h, g$ respectively are assumed atomic rules discharged by an atomic rule. In all cases, the dischargement takes place at the node $f(\mu)$.

\begin{definition}
    $\mathscr{D}$ is \emph{closed} iff all its assumptions are discharged, and it is \emph{open} otherwise.
\end{definition}

\begin{definition}
 Where $\Gamma$ is the set of the undischarged assumptions of $\mathscr{D}$ and $A$ is the conclusion of $\mathscr{D}$, $\mathscr{D}$ is an argument structure \emph{from $\Gamma$ to} (or \emph{for}) \emph{$A$}.
\end{definition}

The notion of \emph{(immediate) sub-structure of $\mathscr{D}$} can be defined as usual. The \emph{substitution of the sub-structure $\mathscr{D}^*$ with the structure $\mathscr{D}^{**}$ in $\mathscr{D}$}---written $\mathscr{D}[\mathscr{D}^{**}/\mathscr{D}^*]$---indicates the argument structure obtained from $\mathscr{D}$ by replacing its sub-structure $\mathscr{D}^*$ with the argument structure $\mathscr{D}^{**}$. Since $\mathscr{D}[\mathscr{D}^{**}/\mathscr{D}^*]$ might not be a sub-structure of $\mathscr{D}$, and since $\mathscr{D}$ is defined as a tree plus discharge functions, when replacing $\mathscr{D}^*$ with $\mathscr{D}^{**}$ in $\mathscr{D}$ one may need to re-define the discharge functions of $\mathscr{D}$, so that assumption formulas or assumed atomic rules discharged at some node $\mu$ in the tree of $\mathscr{D}$ are discharged in $\mathscr{D}[\mathscr{D}^{**}/\mathscr{D}^*]$---if they occur in it---at a node $\mu^*$ which $\mu$ is ‘‘mapped onto" in $\mathscr{D}[\mathscr{D}^{**}/\mathscr{D}]$, whereas some other dischargements in $\mathscr{D}$ might ‘‘disappear" in $\mathscr{D}[\mathscr{D}^{**}/\mathscr{D}]$. I assume that this rough description can be made suitably precise, but I shall abstract from this and only give one example. Take an argument structure

\begin{prooftree}
    \AxiomC{$[\textcolor{red}{A}]$}
    \noLine
    \UnaryInfC{$\mathscr{D}_1$}
    \noLine
    \UnaryInfC{$B$}
    \UnaryInfC{$B \vee C$}
    \AxiomC{$[\textcolor{blue}{B_1}]$}
    \AxiomC{$[\textcolor{blue}{B_2}]$}
    \noLine
    \BinaryInfC{$\mathscr{D}_2$}
    \noLine
    \UnaryInfC{$D$}
    \AxiomC{$[\textcolor{blue}{C}]$}
    \noLine
    \UnaryInfC{$\mathscr{D}_3$}
    \noLine
    \UnaryInfC{$D$}
    \LeftLabel{$\mathscr{D} =$}
    \TrinaryInfC{$\textcolor{green}{D}$}
    \UnaryInfC{$\textcolor{brown}{A \rightarrow D}$}
\end{prooftree}
Here, $B_1$ and $B_2$ are meant to indicate that the assumption $B$ is used twice for deriving $D$ through $\mathscr{D}_2$. Then, the reduced form of $\mathscr{D}$ is $\mathscr{D}[\vee_\rho(\mathscr{D}^*)/\mathscr{D}^*]$, where $\mathscr{D}^*$ is the immediate sub-derivation of $\mathscr{D}$, and $\vee_\rho(\mathscr{D}^*)$ is the reduction of $\mathscr{D}^*$ via the standard reduction for elimination of disjunction, i.e.,

\begin{prooftree}
    \AxiomC{$[\textcolor{red}{A}]$}
    \noLine
    \UnaryInfC{$\mathscr{D}_1$}
    \noLine
    \UnaryInfC{$B_1$}
    \AxiomC{$[\textcolor{red}{A}]$}
    \noLine
    \UnaryInfC{$\mathscr{D}_1$}
    \noLine
    \UnaryInfC{$B_2$}
    \noLine
    \BinaryInfC{$\mathscr{D}_2$}
    \noLine
    \UnaryInfC{$D$}
    \LeftLabel{$\mathscr{D}[\vee_\rho(\mathscr{D}^*)/\mathscr{D}^*] =$}
    \UnaryInfC{$\textcolor{brown}{A \rightarrow D}$}
\end{prooftree}
Thus, the red-to-brown dischargement in $\mathscr{D}$ is re-defined as the double red-to-brown dischargement in $\mathscr{D}[\vee_\rho(\mathscr{D}^*)/\mathscr{D}^*]$, while the blue-to-green dischargements in $\mathscr{D}$ ‘‘disappear" in the dischargements of $\mathscr{D}[\vee_\rho(\mathscr{D}^*)/\mathscr{D}^*]$---see \cite{prawitz1965} for more.

\begin{definition}
    Given $\mathscr{D}$ from $\Gamma = \{\mu_1, ..., \mu_n\}$ to $A$ and $\sigma$ a function from and to argument structures such that $\sigma(\mu_i)$ is a (closed) argument structure with the same conclusion as $\mu_i$ ($i \leq n$), $\mathscr{D}^\sigma = \mathscr{D}[\sigma(\mu_1), ..., \sigma(\mu_n)/\mu_1, ..., \mu_n]$ is called the \emph{(closed) $\sigma$-instance of $\mathscr{D}$}.
\end{definition}

Let me now give an example of argument structure. Level-$0$ rules will be distinguished from assumptions by putting an horizontal bar on top of the former.

\begin{prooftree}
    \AxiomC{$[p \wedge \neg q]_1$}
    \UnaryInfC{$q$}
    \AxiomC{$s$}
    \RightLabel{$[R]_2$}
    \BinaryInfC{$t$}
    \AxiomC{$\overline{[q]}_3$}
    \AxiomC{$\neg \neg p \vee r$}
    \BinaryInfC{$q \rightarrow \neg s$}
    \BinaryInfC{$r$}
    \RightLabel{$3$}
    \UnaryInfC{$t$}
    \AxiomC{}
    \UnaryInfC{$s$}
    \RightLabel{$2$}
    \BinaryInfC{$p$}
    \UnaryInfC{$q \rightarrow (p \vee r)$}
    \RightLabel{$1$}
    \UnaryInfC{$p \wedge \neg q \rightarrow (q \rightarrow (p \vee r))$}
\end{prooftree}
This is open from $s$ and $\neg \neg p \vee r$ to $p \wedge \neg q \rightarrow (q \rightarrow (p \vee r))$. The discharge functions map the non-axiomatic top-node labelled by $p \wedge \neg q$ onto the root-node, the axiomatic top-nodes labelled by $q$ onto the node labelled by $t$ with upper edge labelled by $3$, and the set of edges labelled by $R$ onto the node labelled by $p$ with upper edges labelled by $2$. The structure is an open instance of itself, while a closed instance is the following.

\begin{prooftree}
    \AxiomC{$[p \wedge \neg q]_1$}
    \UnaryInfC{$q$}
    \AxiomC{}
    \UnaryInfC{$s$}
    \RightLabel{$[R]_2$}
    \BinaryInfC{$t$}
    \AxiomC{$\overline{[q]}_3$}
    \AxiomC{}
    \UnaryInfC{$p$}
    \AxiomC{$[q \rightarrow s]_4$}
    \BinaryInfC{$q \vee p$}
    \RightLabel{$4$}
    \UnaryInfC{$\neg \neg p \vee r$}
    \BinaryInfC{$q \rightarrow \neg s$}
    \BinaryInfC{$r$}
    \RightLabel{$3$}
    \UnaryInfC{$t$}
    \AxiomC{}
    \UnaryInfC{$s$}
    \RightLabel{$2$}
    \BinaryInfC{$p$}
    \UnaryInfC{$q \rightarrow (p \vee r)$}
    \RightLabel{$1$}
    \UnaryInfC{$p \wedge \neg q \rightarrow (q \rightarrow (p \vee r))$}
\end{prooftree}
\noindent It should be remarked that, although atomic derivations as per Definition 4 come with embedded dischargements, they can be looked upon as argument structures as per Definition 10 of a special kind: the nodes of the tree are all labelled by atoms, and the discharge-function for assumptions is always the empty function (while the discharge-functions for level-$0$ or higher rules may not be so). In fact, atomic derivations could have been characterised thus, but since below we shall be interested in conditions of derivability in atomic bases, I preferred to define them separately. 

\begin{definition}
    An \emph{inference} is a triple $\langle \langle \mathscr{D}_1, ..., \mathscr{D}_n \rangle, A, \delta \rangle$, where $\delta$ is an extension of the discharge functions associated to the $\mathscr{D}_i$-s ($i \leq n$). The \emph{argument structure associated to the inference}, indicated by the figure
    \begin{prooftree}
        \AxiomC{$\mathscr{D}_1, ..., \mathscr{D}_n$}
        \RightLabel{$\delta$}
        \UnaryInfC{$A$}
    \end{prooftree}
    is obtained by conjoining the trees of $\mathscr{D}_i$-s ($i \leq n$) through a root node $A$, and by adding $\delta$ to the discharges of assumptions of the $\mathscr{D}_i$-s ($i \leq n$). A \emph{rule} is a set of inferences, whose elements are called \emph{instances} of the rule.
\end{definition}

\noindent This definition is inspired by \cite{prawitz1965}. I shall assume that rules amount to recursive sets, which means that they can be described through meta-linguistic schemes. Leaving the set-theoretic notation aside in the case of standard introduction rules in Gentzen's Natural Deduction, the latter are for example given by the following meta-linguistic schemes:
\begin{prooftree}
    \AxiomC{$A$}
    \AxiomC{$B$}
    \RightLabel{($\wedge_I$)}
    \BinaryInfC{$A \wedge B$}
    \AxiomC{$A_i$}
    \RightLabel{($\vee_{I, i}$), $i = 1, 2$}
    \UnaryInfC{$A_1 \vee A_2$}
    \AxiomC{$[A]$}
    \noLine
    \UnaryInfC{$B$}
    \RightLabel{($\rightarrow_I$)}
    \UnaryInfC{$A \rightarrow B$}
    \noLine
    \TrinaryInfC{}
\end{prooftree}
We may consider two additional rules which, sticking this time to Definition 15, are described by the recursive sets

\begin{center}
    $\texttt{Wk} = \{\langle \mathscr{D}, (A \wedge C) \rightarrow B, \emptyset \rangle \ | \ \mathscr{D}$ for $A \rightarrow B\}$
\end{center}
and

\begin{center}
    $\vee_\lambda = \{\langle \mathscr{D}, A, \emptyset \rangle \ | \ \mathscr{D}$ for $A \vee B\}$.
\end{center}
In more familiar Natural Deduction notation, these correspond to the meta-linguistic schemes

\begin{prooftree}
    \AxiomC{$A \rightarrow B$}
    \RightLabel{$\texttt{Wk}$}
    \UnaryInfC{$(A \wedge C) \rightarrow B$}
    \AxiomC{$A \vee B$}
    \RightLabel{$\vee_\lambda$}
    \UnaryInfC{$A$}
    \noLine
    \BinaryInfC{}
\end{prooftree}

\begin{definition}
    $\mathscr{D}$ is \emph{canonical} iff it is associated to an instance of an introduction rule. It is \emph{non-canonical} otherwise.
\end{definition}

\begin{definition}
    Given a rule $R$, a \emph{reduction for $R$} is a mapping $\phi$ from and to argument structures such that $\phi$ is defined on some sub-set $\mathbb{D}$ of the argument structures associated to instances of $R$ and, $\forall \mathscr{D} \in \mathbb{D}$,
    \begin{itemize}
        \item[(a)] $\mathscr{D}$ is from $\Gamma$ to $A \Longrightarrow \phi(\mathscr{D})$ is from some $\Gamma^* \subseteq \Gamma$ to $A$;
        \item[(b)] $\forall \sigma$, $\phi$ is defined on $\mathscr{D}^\sigma$ and $\phi(\mathscr{D}^\sigma) = \phi(\mathscr{D})^\sigma$.
    \end{itemize}
\end{definition}

\noindent I shall use the following notation: $\mathfrak{J}$ indicates a set of reductions, whereas $\mathfrak{J}^+$ indicates an extension of $\mathfrak{J}$, i.e., $\mathfrak{J} \subseteq \mathfrak{J}^+$. One can easily see that the standard reductions for elimination of conjunction, disjunction and implication, as used in proofs of normalisation for Natural Deduction \cite{prawitz1965}, all cope with the requirements of Definition 17. We could associate the non-introduction rules $\texttt{Wk}$ and $\vee_\lambda$ defined above to the respective reductions

\begin{prooftree}
    \AxiomC{$[A]_1$}
    \noLine
    \UnaryInfC{$\mathscr{D}$}
    \noLine
    \UnaryInfC{$B$}
    \RightLabel{$1$}
    \UnaryInfC{$A \rightarrow B$}
    \RightLabel{$\texttt{Wk}$}
    \UnaryInfC{$(A \wedge C) \rightarrow B$}
    \AxiomC{$\stackrel{\texttt{Wk}_\rho}{\Longrightarrow}$}
    \noLine
    \UnaryInfC{}
    \AxiomC{$[A \wedge C]_1$}
    \UnaryInfC{$A$}
    \noLine
    \UnaryInfC{$\mathscr{D}$}
    \noLine
    \UnaryInfC{$B$}
    \RightLabel{$1$}
    \UnaryInfC{$(A \wedge C) \rightarrow B$}
    \noLine
    \TrinaryInfC{}
\end{prooftree}
and

\begin{prooftree}
    \AxiomC{$\mathscr{D}$}
    \noLine
    \UnaryInfC{$A$}
    \UnaryInfC{$A \vee B$}
    \RightLabel{$\vee_\lambda$}
    \UnaryInfC{$A$}
    \AxiomC{$\stackrel{\vee_{\lambda, \ \rho}}{\Longrightarrow}$}
    \noLine
    \UnaryInfC{}
    \AxiomC{$\mathscr{D}$}
    \noLine
    \UnaryInfC{$A$}
    \noLine
    \TrinaryInfC{}
\end{prooftree}
\noindent Observe that both the standard reductions for conjunction, disjunction and implication, and $\texttt{Wk}_\rho$ and $\vee_{\lambda, \ \rho}$, are defined on \emph{proper} sub-sets of the set of the argument structures associated to all the instances of the corresponding rules---i.e., those where the argument structure for the major premise is in canonical form (for only one of the disjuncts in the case of $\vee_{\lambda, \ \rho}$).

\begin{definition}
    $\mathscr{D}$ \emph{immediately reduces to $\mathscr{D}^*$ relative to $\mathfrak{J}$} iff, for some $\mathscr{D}^{**}$ sub-structure of $\mathscr{D}$, there is $\phi \in \mathfrak{J}$ such that $\mathscr{D}^* = \mathscr{D}[\phi(\mathscr{D}^{**})/\mathscr{D}^{**}]$. The relation of $\mathscr{D}$ \emph{reducing to $\mathscr{D}^*$ relative to $\mathfrak{J}$} is the reflexive-transitive closure of the relation of immediate reducibility of $\mathscr{D}$ to $\mathscr{D}^*$.
\end{definition}

\noindent I will use the following notation: $\mathscr{D} \leq_\mathfrak{J} \mathscr{D}^*$ means that $\mathscr{D}$ reduces to $\mathscr{D}^*$ relative to $\mathfrak{J}$.

\begin{definition}
    An \emph{argument} is a pair $\langle \mathscr{D}, \mathfrak{J} \rangle$.
\end{definition}

\begin{definition}
    $\langle \mathscr{D}, \mathfrak{J} \rangle$ is \emph{$n$-valid on $\mathfrak{B}$} iff $\mathfrak{B} \in \mathbb{B}^n$ and
    \begin{itemize}
        \item $\mathscr{D}$ is closed $\Longrightarrow$
        \begin{itemize}
            \item the conclusion of $\mathscr{D}$ is atomic $\Longrightarrow \mathscr{D} \leq_\mathfrak{J} \mathscr{D}^*$ with $\mathscr{D}^* \in \texttt{DER}_{\mathfrak{B}}$ closed;
            \item the conclusion of $\mathscr{D}$ is not atomic $\Longrightarrow \mathscr{D} \leq_\mathfrak{J} \mathscr{D}^*$ for $\mathscr{D}^*$ canonical with immediate sub-structures $n$-valid on $\mathfrak{B}$ when paired with $\mathfrak{J}$;
        \end{itemize}
        \item $\mathscr{D}$ is open $\Longrightarrow \forall \sigma, \forall \mathfrak{J}^+ \supseteq \mathfrak{J}$ and $\forall \mathfrak{C} \supseteq_n \mathfrak{B}$, if, $\forall \mu \in \Gamma$, $\langle \sigma(\mu), \mathfrak{J}^+ \rangle$ is $n$-valid on $\mathfrak{C}$, then $\langle \mathscr{D}^\sigma, \mathfrak{J}^+ \rangle$ is $n$-valid on $\mathfrak{C}$.
    \end{itemize}
\end{definition}

\noindent I will use the following notation: $\Gamma \models^\alpha_{\mathfrak{B}, \ n} A$ means that there is $\langle \mathscr{D}, \mathfrak{J} \rangle$ from $\Gamma$ to $A$ which is $n$-valid on $\mathfrak{B}$. Based on the reduction $\vee_{\lambda, \ \rho}$ for the rule $\vee_\lambda$ above, it is easy to see that, for every $n$, $A \vee B \models^\alpha_{\mathfrak{B}, \ n} A$ as soon as $\models^\alpha_{\mathfrak{B}, \ n} \neg B$.

\begin{proposition}[Equivalence result]
    For every $n$ and every $\mathfrak{B} \in \mathbb{B}^n$:
    \begin{itemize}
        \item[(a)] $\forall A \in \texttt{ATOM}_\mathscr{L} \ (\models^\alpha_{\mathfrak{B}, \ n} A \Longleftrightarrow \ \vdash_\mathfrak{B} A)$;
        \item[(b)] $\models^\alpha_{\mathfrak{B}, \ n} \bot \Longleftrightarrow \ \forall A \in \texttt{ATOM}_\mathscr{L} \ (\vdash_\mathfrak{B} A)$;
        \item[(c)] $\models^\alpha_{\mathfrak{B}, \ n} A \wedge B \Longleftrightarrow \ \models^\alpha_{\mathfrak{B}, \ n} A$ and $\models^\alpha_{\mathfrak{B}, \ n} B$;
        \item[(d)] $\models^\alpha_{\mathfrak{B}, \ n} A \vee B \Longleftrightarrow \ \models^\alpha_{\mathfrak{B}, \ n} A$ or $\models^\alpha_{\mathfrak{B}, \ n} B$;
        \item[(e)] $\models^\alpha_{\mathfrak{B}, \ n} A \rightarrow B \Longleftrightarrow \ A \models^\alpha_{\mathfrak{B}, \ n} B$;
        \item[(f)] $\Gamma \models^\alpha_{\mathfrak{B}, \ n} A \Longleftrightarrow \forall \mathfrak{C} \supseteq_n \mathfrak{B} \ (\Gamma \models^\alpha_{\mathfrak{C}, \ n} A)$---i.e. the notion is monotonic;
        \item[(g)] \emph{(Admissibility clause)} $\Gamma \models^\alpha_{\mathfrak{B}, \ n} A \Longleftrightarrow \ \forall \mathfrak{C} \supseteq_n \mathfrak{B} \ (\models^\alpha_{\mathfrak{C}, \ n} \Gamma \Longrightarrow \ \models^\alpha_{\mathfrak{C}, \ n} A)$.
    \end{itemize}
\end{proposition}

\begin{proof}
    The only non-trivial case is direction ($\Longleftarrow$) of (g).\footnote{A similar proof is to be found in \cite[Lemma 4.11]{stafford1}. Since Stafford's discussion has no restrictions on $\Gamma$ being finite, in that context the result amounts to proving that Prawitz's reducibility semantics enjoys compact monotonicity. Another similar proof is found in \cite[pp. 153-154]{piccolomininote}, although there the proof is referred to non-monotonic reducibility semantics, and used for showing that, when suitable conditions are met, classical logic can be justified over Prawitz's reducibility semantics.} Assume $\forall \mathfrak{C} \supseteq_n \mathfrak{B} \ (\models^{\alpha}_{\mathfrak{C}, \ n} \Gamma \Longrightarrow \ \models^\alpha_{\mathfrak{C}, \ n} A)$. Let $\Gamma = \{A_1, ..., A_m\}$, let $\mathbb{F}$ be the set of reductions according to Definition 17, and let $\mathscr{D}$ be the argument structure

    \begin{prooftree}
        \AxiomC{$A_1$}
        \AxiomC{$\dots$}
        \AxiomC{$A_m$}
        \TrinaryInfC{$A$}
    \end{prooftree}
    We are going to establish that $\Gamma \models^\alpha_{\mathfrak{C}, \ n} A$ by showing that $\langle \mathscr{D}, \mathbb{F} \rangle$ is $n$-valid on $\mathfrak{B}$. Take any $\sigma$, any $\mathbb{F}^+$, and any $\mathfrak{C} \supseteq_n \mathfrak{B}$ such that $\langle \sigma(A_i), \mathbb{F}^+ \rangle \ (i \leq m)$ is $n$-valid on $\mathfrak{C}$. We must prove that $\langle \mathscr{D}^\sigma, \mathbb{F}^+ \rangle$ is $n$-valid on $\mathfrak{C}$, i.e., that

    \begin{prooftree}
        \AxiomC{$\sigma(A_1)$}
        \AxiomC{$\dots$}
        \AxiomC{$\sigma(A_m)$}
        \TrinaryInfC{$A$}
    \end{prooftree}
    is $n$-valid on $\mathfrak{C}$ when paired with $\mathbb{F}^+$. We observe first of all that, since $\mathbb{F}$ is the set of all the reductions, we have $\mathbb{F}^+ = \mathbb{F}$, so $\langle \sigma(A_i), \mathbb{F} \rangle \ (i \leq m)$ is $n$-valid on $\mathfrak{C}$. Also, by our initial assumption we have $\models^\alpha_{\mathfrak{C}, \ n} A$, i.e., there is closed $\mathscr{D}_A$ for $A$ such that, for some $\mathfrak{J}$, $\langle \mathscr{D}_A, \mathfrak{J} \rangle$ is $n$-valid on $\mathfrak{C}$. Since validity of argument structures is monotonic over extensions of sets of reductions, we then have that $\langle \mathscr{D}_A, \mathbb{F} \rangle$ is $n$-valid on $\mathfrak{C}$. Now, let $\phi^\sigma$ be the mapping that is defined for $\mathscr{D}^\sigma$ only, and maps it to $\mathscr{D}_A$. Then $\phi^\sigma$ is a reduction (for the rule whose only instance is $\mathscr{D}^\sigma$), i.e., it satisfies conditions (a) and (b) of Definition 17: condition (a) holds trivially while, for condition (b), since $\mathscr{D}^\sigma$ and $\mathscr{D}_A$ are closed, for every mapping $\tau$ from formulas to argument structures as required by Definition 14, we have that $\phi^\sigma((\mathscr{D}^\sigma)^\tau) = \phi^\sigma(\mathscr{D}^\sigma) = \mathscr{D}_A = \mathscr{D}^\tau_A = \phi^\sigma(\mathscr{D}^\sigma)^\tau$. So $\phi^\sigma \in \mathbb{F}$, and $\mathscr{D}^\sigma \leq_\mathbb{F} \mathscr{D}_A$, i.e., $\mathscr{D}^\sigma$ reduces modulo $\mathbb{F}$ to a closed argument structure for $A$ which is $n$-valid over $\mathfrak{C}$ when paired with $\mathbb{F}$, which means that $\langle \mathscr{D}^\sigma, \mathbb{F} \rangle$ is $n$-valid over $\mathfrak{C}$ too, as desired.\footnote{A previous rendering of this proof was much more baroque and difficult to read. I am grateful to one of the reviewers for proposing an alternative and much smoother rendering, which I was glad to implement in my paper.}
\end{proof}

Below I shall refer many times to Proposition 5---especially to its point (g). This is why, to help readability, I gave them a name---Proposition 5 is the \emph{equivalence result}, while point (g) is the \emph{admissibility clause}.

Concerning the proof of the admissibility clause, observe that, in principle, there might be a different $\phi^\sigma$ for each different sequence $\sigma(A_1), ..., \sigma(A_m)$ such that $\langle \sigma(A)_i, \mathbb{F} \rangle$ is closed valid for $A_i$ on $\mathfrak{C} \ (i \leq m)$, since there might in principle be a different $\mathscr{D}_A$ such that $\langle \mathscr{D}_A, \mathfrak{J} \rangle$ is closed valid for $A$ on $\mathfrak{C}$ associated to each such sequence. With classical logic in the meta-language, however, the situation is much smoother: either there is $A_i$ such that there is no closed argument for $A_i$ valid on $\mathfrak{C}$, in which case $\phi^\sigma$ is the empty function, or for every $A_i$ there is a closed argument for $A_i$ valid on $\mathfrak{C}$, in which case there will surely be at least one closed argument for $A$ valid on $\mathfrak{C}$, so $\phi^\sigma$ can be the constant function.\footnote{One may say that the reduction $\mathbb{F}$ in the proof of the admissibility clause is not ‘‘constructive enough". More constructive examples of reductions \emph{can} be given---these are essentially adaptations from the incompleteness results proved in \cite{piechadecampossanz, piechaschroeder-heisterdecampossanz, piechaschroeder-heister2019}. The discussion of this topic would have led me too far away, so I may want to deal with it in future works---for a partial treatment, see \cite{piccolomininote}. Concerning the issue of what a ‘‘good" reduction is, the reader may refer to \cite{ayhangoodreductions}. Let me also point out that, as done in \cite{schroeder-heister2006}, reductions could be also defined in terms of (sequences of) pairs of argument structures, i.e., sequence of pairs $\langle \mathscr{D}, \mathscr{D}^* \rangle$ where the first element of the $i + 1$-th pair has at most the assumptions and the same conclusion as the second element in $i$-th pair. One final observation is that, besides few very quick remarks, I shall not deal here with ‘‘deviant" reducibility semantics where potential readings of what a reduction should be, stricter than those mentioned above, make the admissibility clause fail.}

\begin{definition}
    $\langle \mathscr{D}, \mathfrak{J} \rangle$ is \emph{$n$-valid} iff, $\forall \mathfrak{B} \in \mathbb{B}^n$, $\langle \mathscr{D}, \mathfrak{J} \rangle$ is $n$-valid on $\mathfrak{B}$.
\end{definition}

\noindent I will use the following notation: $\Gamma \models^\alpha_n A$ means that there is an $n$-valid $\langle \mathscr{D}, \mathfrak{J} \rangle$ from $\Gamma$ to $A$. Based on the reduction $\texttt{Wk}_\rho$ for the rule $\texttt{Wk}$ above, it is easy to see that, for every $n$, $A \rightarrow B \models^\alpha_n (A \wedge C) \rightarrow B$.

\begin{proposition}[Schroeder-Heister, \cite{schroeder-heister2006}]

$\Gamma \models^\alpha_n A \Longleftrightarrow \Gamma \models^\alpha_{\mathfrak{B}^\emptyset, \ n} A$.
\end{proposition}

\begin{corollary}
    The following facts hold:
    \begin{itemize}
    \item[(a)] $\Gamma \models^\alpha_n A \Longleftrightarrow \forall \mathfrak{B} \in \mathbb{B}^n \ (\Gamma \models^\alpha_{\mathfrak{B}, \ n} A)$;
    \item[(b)] $\forall \mathfrak{B} \in \mathbb{B}^n \ (\Gamma \models^\alpha_{\mathfrak{B}, \ n} A) \Longleftrightarrow \forall \mathfrak{B} \in \mathbb{B}^n \ (\models^\alpha_{\mathfrak{B}, \ n} \Gamma \Longrightarrow \ \models^\alpha_{\mathfrak{B}, \ n} A)$.
    \end{itemize}
\end{corollary}

\begin{proof}
   As for (a), ($\Longrightarrow$) is trivial. For ($\Longleftarrow$), if $\forall \mathfrak{B} \in \mathbb{B}^n \ ( \Gamma \models^\alpha_{\mathfrak{B}, \ n} A)$ then, in particular, $\Gamma \models^\alpha_{\mathfrak{B}^\emptyset, \ n} A$ and, by Proposition 6, $\Gamma \models^\alpha_n A$. As for (b), ($\Longrightarrow$) if $\forall \mathfrak{B} \in \mathbb{B}^n \ ( \Gamma \models^\alpha_{\mathfrak{B}, \ n} A)$, then $\Gamma \models^\alpha_{\mathfrak{B}^\emptyset, \ n} A$. By the admissibility clause plus Proposition 1, this means $\forall \mathfrak{B} \in \mathbb{B}^n \ (\models^\alpha_{\mathfrak{B}, \ n} \Gamma \Longrightarrow \ \models^\alpha_{\mathfrak{B}, \ n} A)$. ($\Longleftarrow$) Instantiate $\forall \mathfrak{B} \in \mathbb{B}^n \ (\models_{\mathfrak{B}, \ n} \Gamma \Longrightarrow \ \models_{\mathfrak{B}, \ n} A)$ on any $\mathfrak{B} \in \mathbb{B}^n$. By the admissibility clause, $\models_{\mathfrak{B}, \ n} \Gamma \Longrightarrow \ \models_{\mathfrak{B}, \ n} A$ implies $\Gamma \models_{\mathfrak{B}, \ n} A$. We can now universally quantify over $\mathfrak{B} \in \mathbb{B}^n$. 
\end{proof}

\section{Comparison of the three approaches}

In what follows, I shall be comparing the three approaches to proof-theoretic semantics mentioned so far.

\subsection{Reducibility semantics and standard base semantics}

Definition 7 and the equivalence result have the obvious effect of making standard base semantics and reducibility semantics equivalent, both at the level of consequence over a base, and at the level of logical consequence.

\begin{theorem}
    $\Gamma \models_{\mathfrak{B}, \ n} A \Longleftrightarrow \Gamma \models^\alpha_{\mathfrak{B}, \ n} A$.
\end{theorem}

\begin{theorem}
    $\Gamma \models_n A \Longleftrightarrow \Gamma \models^\alpha_n A$.
\end{theorem}

\noindent So, we can transfer to reducibility semantics all the (in)completeness results which have been proved for base semantics in the standard reading. Below, we shall discuss some incompleteness results for intuitionistic logic $\texttt{IL}$. For the moment the reader may refer to \cite{piechaschroeder-heisterdecampossanz, piechaschroeder-heister2019, stafford1}, and rely on the following result.

\begin{theorem}
    $\exists \Gamma \ \exists A \ (\Gamma \models^\alpha_n A$ and $\Gamma \not\vdash_{\emph{\texttt{IL}}} A)$.
\end{theorem}

\noindent Despite their simplicity, Theorem 1 and Theorem 2 provide an interesting connection between Prawitz's reducibility semantics and base semantics in the standard reading. For, the latter can be looked at as ‘‘extracted" from the former by dropping argument structures and reductions out, via a consequence relation defined over (sets of) formulas outright, rather than in a derivative way as existence of suitable valid arguments. Conversely, Prawitz's reducibility semantics can be understood as obtained from base semantics in the standard reading, by decorating the formula-based consequence relation with argument structures and reductions which ‘‘witness" that such a relation holds. Hence, one may naturally wonder whether, whenever $\Gamma \models_{\mathfrak{B}, \ n} A$ or $\Gamma \models_n A$ hold, a ‘‘witness" of this can be found so that $\Gamma \models^\alpha_{\mathfrak{B}, \ n} A$ or $\Gamma \models^\alpha_n A$ hold too, and conversely whether, whenever $\Gamma \models^\alpha_{\mathfrak{B}, \ n} A$ or $\Gamma \models^\alpha_n A$ hold, the ‘‘witness" can be safely removed without loss of constructivity.

As I see it, the point raised here is best appreciated when formulated in a style similar to the one at play in another major constructivist approach, i.e., Martin-L\"{o}f's (intuitionistic) type theory \cite{martin-loef}. Here, following the formulas-as-types conception and the Curry-Howard correspondence inspiring it \cite{howard}, we start with forms of \emph{judgement}

\begin{center}
    $a : A$
\end{center}
where $a$ is a proof-object and $A$ is a type/proposition, so what the judgement expresses is that the proof-object $a$ is of type $A$, i.e., it is a proof-object of proposition $A$ and, conversely, that proposition $A$, understood as a type, i.e., as a class of proof-objects, is inhabited. Based on this, one can introduce forms of judgement

\begin{center}
    $A \ \texttt{true}$
\end{center}
where $\texttt{true}$ is explained by the rule

\begin{prooftree}
    \AxiomC{$a : A$}
    \RightLabel{$\texttt{T}$}
    \UnaryInfC{$A \ \texttt{true}$}
\end{prooftree}
Now, Prawitz's reducibility semantics runs in very much the same way, in that $\models^\alpha_\mathfrak{B}$ and $\models^\alpha$ are always explained in terms of (existence of) proof-objects witnessing them, similarly to what happens in Martin-L\"{o}f's intuitionistic type theory, when explaining $\texttt{true}$ through the rule $\texttt{T}$. In standard base semantics, instead, $\models_\mathfrak{B}$ and $\models$ are defined at the sentence level directly so that, to put it in type-theoretic terms, it is as if we started with $\texttt{true}$ outright, without requiring it to be explained by $\texttt{T}$. To distinguish this picture from the previous one, we may write $\texttt{true}^*$ instead of $\texttt{true}$. Now, the question I raised above---whether from $\Gamma \models_\mathfrak{B}$ or $\Gamma \models A$ one can go to $\Gamma \models^\alpha_\mathfrak{B} A$ or $\Gamma \models^\alpha A$ and vice versa---would become the question whether one of the following rules is valid---leaving atomic bases out and using a sequent-style notation:

\begin{prooftree}
    \AxiomC{$\Gamma \ \textbf{true} \Rightarrow A \ \texttt{true}$}
    \UnaryInfC{$\Gamma \ \textbf{true}^* \Rightarrow A \ \texttt{true}^*$}
    \AxiomC{$\Gamma \ \textbf{true}^* \Rightarrow A \ \texttt{true}^*$}
    \UnaryInfC{$\Gamma \ \textbf{true} \Rightarrow A \ \texttt{true}$}
    \noLine
    \BinaryInfC{}
\end{prooftree}
---where $\textbf{true}^{(*)}$ means that all the elements of $\Gamma$ are $\texttt{true}^{(*)}$---namely, whether one of the following is valid:

\begin{prooftree}
    \AxiomC{$\textbf{x} : \Gamma \Rightarrow b(\textbf{x}) : A(\textbf{x})$}
    \UnaryInfC{$\Gamma \ \textbf{true}^* \Rightarrow A \ \texttt{true}^*$}
    \AxiomC{$\Gamma \ \textbf{true}^* \Rightarrow A \ \texttt{true}^*$}
    \UnaryInfC{$\textbf{x} : \Gamma \Rightarrow b(\textbf{x}) : A(\textbf{x})$}
    \noLine
    \BinaryInfC{}
\end{prooftree}
where $\textbf{x} : \Gamma \Rightarrow b(\textbf{x}) : A(\textbf{x})$ means $b(\textbf{a}) : A(\textbf{a})$ whenever $\textbf{a}$ is a sequence of proof-objects for the elements of $\Gamma$. Theorem 1 and Theorem 2 answer positively to these questions. At the same time, they imply that, under the given limitations (mostly concerning the restriction to finite $\Gamma$-s), Prawitz's reducibility semantics and standard base semantics are ‘‘structurally" identical, meaning that the former enjoys \emph{mutatis mutandis} the same properties as those employed in \cite{piechaschroeder-heisterdecampossanz, piechaschroeder-heister2019} to prove incompleteness of $\texttt{IL}$ with respect to the latter---properties which I shall come back to below. So, as stated by Theorem 3, $\texttt{IL}$ is also incomplete over Prawitz's reducibility semantics in the version at issue here---essentially the same line of thought is used in \cite{piccolominibook} to prove incompleteness of $\texttt{IL}$ over Prawitz's Theory of Grounds \cite{prawitz2015}.

These full-equivalence results are mainly due to the fact that both Prawitz's reducibility semantics and base semantics in the standard reading are ‘‘introduction-based", i.e., they both ultimately rely on the idea that the meaning of a logical constant $\kappa$ is given by the \emph{conditions} for formulas having $\kappa$ as main sign to hold. This is the reason why Theorem 1 can be proved by easy induction on the complexity of formulas in the closed case, and then by ‘‘closure" in the open case, using the admissibility clause---then, Theorem 2 follows from the equivalence result, and Theorem 3 follows straightforwardly from Theorem 2 plus the incompleteness results proved in \cite{piechaschroeder-heisterdecampossanz, piechaschroeder-heister2019}.

\subsection{Reducibility semantics and Sandqvist's base semantics}

However, this is also the reason why such a smooth inductive reasoning, hence full-equivalence, do not apply when comparing Prawitz's reducibility semantics and base semantics \emph{à la} Sandqvist. The reasoning breaks down since Sandqvist explains $\vee$ in an ‘‘elimination-based" way. For, clause (c) in Definition 8 is in fact nothing but a (monotonic) semantic rendering of the standard elimination rule for $\vee$, restricted to atomic minor premises. But the question may arise also here whether, as in the previous case, whenever $\Gamma \models^s_{\mathfrak{B}, \ n} A$ or $\Gamma \models^s_n A$ hold, a ‘‘witness" of this can be found so that $\Gamma \models^\alpha_{\mathfrak{B}, \ n} A$ or $\Gamma \models^\alpha_n A$ hold too, and vice versa, whether from $\Gamma \models^\alpha_{\mathfrak{B}, \ n} A$ or $\Gamma \models^\alpha_n A$ one can conclude $\Gamma \models^s_{\mathfrak{B}, \ n} A$ or $\Gamma \models^s_n A$. Because smooth induction and full-equivalence are now ruled out, anyway, the most one can expect in this context is establishing \emph{to what extent} the back-and-forth is possible. In turn, this can be understood as establishing the extent to which the smooth inductive reasoning used for proving Theorem 1 can be applied in the case of a comparison between Prawitz's reducibility semantics and Sandqvist's base semantics.

Going inductively from Prawitz's reducibility semantics to base semantics in Sandqvist's variant might seem, at first, not beyond reach. Omitting details to be found in the proof of Theorem 4, suppose $A \vee B$ holds in Prawitz's sense on some base. Then, for every extension of the base, either $A$ or $B$ hold on that extension. Assume that atom $C$ is consequence of both $A$ and $B$ in Sandqvist's sense on that extension. If we can grant inductively that the holding of $A$ and $B$ in Prawitz's sense on the extension implies the holding of $A$ and $B$ in Sandqvist's sense on the extension (which is however something that we \emph{cannot} grant, as we shall see in a while), we can conclude $C$ holds in Sandqvist's sense on the extension too which, by arbitrariness of the extension, means that $A \vee B$ holds in Sandqvist's sense in the original base.

The inverse, however, seems not to hold: from the fact that, for every atom $C$ which is consequence of $A$ and $B$ in Sandqvist's sense on some extension of a given base, $C$ holds in Sandqvist's sense on the extension, we cannot conclude that $A$ or $B$ hold on the base in Sandqvist's sense, so we cannot apply any potential induction for going from Sandqvist's base semantics to Prawitz's reducibility semantics---which is, incidentally, one of the selling points of Sandqvist's approach as concerns completeness of intuitionistic logic, see \cite{sandqvist}.

It is thus clear that an inductive full-equivalence proof is broken also for the other propositional constants. This does not exclude, however, that \emph{assuming} the condition that we can go from Sandqvist to Prawitz everywhere, we can keep conditionally, for \emph{all} propositional constants, the possibility of going the other way around everywhere as well.

The aforementioned condition is required mainly because of how the clause for $\rightarrow$ and that for $\Gamma \neq \emptyset$ work. Again omitting details to be found in the proof of Theorem 4, and limiting ourselves to the implicational case, suppose that $A \rightarrow B$ holds in Prawitz's sense on some base. Therefore, for any extension of the base, if $A$ holds in Prawitz's sense on the extension, also $B$ holds in Prawitz's sense on the extension. We must prove that the same holds on the extension also for Sandqvist's variant, namely that, whenever $A$ holds in Sandqvist's sense on the extension, then also $B$ holds in Sandqvist's sense on the extension. This obtains in turn when we assume that, for every extension and every $\Gamma$ and $A$, we can go from $A$ being consequence of $\Gamma$ in Sandqvist's sense on the extension to $A$ being consequence of $\Gamma$ in Prawitz's sense on the extension.\footnote{This is a stronger claim than what is actually needed for proving Theorem 4. We could limit ourselves to assuming that the implication holds only for valid formulas, rather than for consequences (i.e., with $\Gamma = \emptyset$). However, the stronger claim clearly implies the weaker, and is needed in the proof of Theorem 13, which is why I favoured it.}

Thus, Prawitz's reducibility semantics and base semantics in Sandqvist's reading \emph{are} connected, but not ‘‘point-wise", i.e., assigning to \emph{each} instance $\Gamma \models^s_{\mathfrak{B}, \ n} A$ of Sandqvist's consequence notion the corresponding instance $\Gamma \models^\alpha_{\mathfrak{B}, \ n} A$ of Prawitz's consequence notion, and vice versa. The connection is, so to say, ‘‘conditionally global". As said, this certainly depends on the structural difference given by Prawitz's ‘‘introduction-based", and Sandqvist's ‘‘elimination-based" explanations of $\vee$. This ‘‘conditionally global" connection is however not without value, as we shall see. Before turning to that, however, let me prove the results which I have been illustrating informally so far. I will use the following notation: $A^* \preceq A$ indicates that the logical complexity of $A^*$ is smaller than or equal to that of $A$. Consider now:

\begin{equation}
    \forall \Gamma \ \forall A \ \forall \mathfrak{C} \supseteq_n \mathfrak{B} \ (\Gamma \models^s_{\mathfrak{C}, \ n} A \Longrightarrow \Gamma \models^\alpha_{\mathfrak{C}, \ n} A).
\end{equation}

\begin{theorem}
    $(1) \Longrightarrow \forall \Gamma \ \forall A \ \forall \mathfrak{C} \supseteq_n \mathfrak{B} \ (\Gamma \models^\alpha_{\mathfrak{C}, \ n} A \Longrightarrow \Gamma \models^s_{\mathfrak{C}, \ n} A)$.
\end{theorem}

\begin{proof}
Assume (1), and suppose $\Gamma = \emptyset$. We proceed by induction on complexity of $A$:

    \begin{itemize}
        \item take any arbitrary $\mathfrak{C} \supseteq_n \mathfrak{B}$ and suppose $A \in \texttt{ATOM}$. Then we have
        \begin{center}
        $\models^\alpha_{\mathfrak{C}, \ n} A$ iff $\vdash_\mathfrak{C} A$ iff $\models^s_{\mathfrak{C}, \ n} A$.
        \end{center}
        Hence, \emph{a fortiori},
        \begin{center}
        $\forall \mathfrak{C} \supseteq_n \mathfrak{B} \ (\models^\alpha_{\mathfrak{C}, \ n} A \Longrightarrow \ \models^s_{\mathfrak{C}, \ n} A)$;
        \end{center}
        Observe that we do not have to worry about $A = \bot$ via Proposition 2;
        \item assume the induction hypothesis
        \begin{center}
            $\forall A^* \prec A \ \forall \mathfrak{C} \supseteq_n \mathfrak{B} \ (\models^\alpha_{\mathfrak{C}, \ n} A^* \Longrightarrow \ \models^s_{\mathfrak{C}, \ n} A^*)$.
        \end{center}
        The case for $\wedge$ is trivial. As for the others:
    \begin{itemize}
        \item[($\vee$)] take any arbitrary $\mathfrak{C} \supseteq_n \mathfrak{B}$ and suppose $A = B \vee C$. So, by points (d) and (f) of the equivalence result
        \begin{center}
        $\models^\alpha_{\mathfrak{C}, \ n} B \vee C$ iff $\forall \mathfrak{C}^* \supseteq_n \mathfrak{C} \ (\models^\alpha_{\mathfrak{C}^*, \ n} B \vee C)$ iff $\forall \mathfrak{C}^* \supseteq_n \mathfrak{C} \ (\models^\alpha_{\mathfrak{C}^*, \ n} B$ or $\models^\alpha_{\mathfrak{C}^*, \ n} C)$.
        \end{center}
        For any arbitrary $\mathfrak{C}^* \supseteq_n \mathfrak{C}$ and any arbitrary $D \in \texttt{ATOM}_{\mathscr{L}}$, assume $B \models^s_{\mathfrak{C}^*, \ n} D$ and $C \models^s_{\mathfrak{C}^*, \ n} D$. Instantiate on $\mathfrak{C}^*$ the previous
        \begin{center}
            $\forall \mathfrak{C}^* \supseteq \mathfrak{C} \ (\models^\alpha_{\mathfrak{C}^*, \ n} B$ or $\models^\alpha_{\mathfrak{C}^*, \ n} C)$
        \end{center}
        and assume $\models^\alpha_{\mathfrak{C}^*, \ n} B$. Instantiate on $B$ the induction hypothesis, so to obtain $\models^s_{\mathfrak{C}^*, \ n} B$ which, with $B \models^s_{\mathfrak{C}^*, \ n} D$, yields $\models^s_{\mathfrak{C}^*, \ n} D$. The same can be obtained by assuming $\models^\alpha_{\mathfrak{C}^*, \ n} C$. Discharging the assumptions $\models^\alpha_{\mathfrak{C}^*, \ n} B$ and $\models^\alpha_{\mathfrak{C}^*, \ n} C$, we obtain $\models^s_{\mathfrak{C}^*, \ n} D$ and, discharging the assumptions $B \models^s_{\mathfrak{C}^*, \ n} D$ and $C \models^s_{\mathfrak{C}^*, \ n} D$, we have
        \begin{center}
            $B \models^s_{\mathfrak{C}^*, \ n} D$ and $C \models^s_{\mathfrak{C}^*, \ n} D \Longrightarrow \ \models^s_{\mathfrak{C}^*, \ n} D$.
        \end{center}
        $\mathfrak{C}^*$ and $D$ do not occur free in any undischarged assumption so we obtain
        \begin{center}
            $\forall \mathfrak{C}^* \supseteq_n \mathfrak{C} \ \forall D \in \texttt{ATOM}_{\mathscr{L}} \ (B \models^s_{\mathfrak{C}^*, \ n} D$ and $C \models^s_{\mathfrak{C}^*, \ n} D \Longrightarrow \ \models^s_{\mathfrak{C}^*, \ n} D)$
        \end{center}
        which, by point (c) in Definition 8, means $\models^s_{\mathfrak{C}, \ n} B \vee C$. So we have
        \begin{center}
        $\forall \mathfrak{C} \supseteq_n \mathfrak{B} \ (\models^\alpha_{\mathfrak{C}, \ n} B \vee C \Longrightarrow \ \models^s_{\mathfrak{C}, \ n} B \vee C)$;
        \end{center}
        \item[($\rightarrow$)] take any arbitrary $\mathfrak{C} \supseteq_n \mathfrak{B}$ and suppose $A = B \rightarrow C$. Then, by point (e) and (g) of the equivalence result,
            \begin{center} 
            $\models^\alpha_{\mathfrak{C}, \ n} B \rightarrow C$ iff $B \models^\alpha_{\mathfrak{C}, \ n} C$ iff $\forall \mathfrak{C}^* \supseteq_n \mathfrak{C} \ (\models^\alpha_{\mathfrak{C}^*, \ n} B \Longrightarrow \ \models^\alpha_{\mathfrak{C}^*, \ n} C)$.
            \end{center} 
            Take now any arbitrary $\mathfrak{C}^* \supseteq_n \mathfrak{C}$ and assume $\models^s_{\mathfrak{C}^*, \ n} B$. Instantiating (1) on $\emptyset$, $B$ and $\mathfrak{C}^* \supseteq_n \mathfrak{C} \supseteq_n \mathfrak{B}$, we get
            \begin{center} 
            $\models^s_{\mathfrak{C}^*, \ n} B \Longrightarrow \ \models^\alpha_{\mathfrak{C}^*, \ n} B$
            \end{center} 
            which yields $\models^\alpha_{\mathfrak{C}^*, \ n} B$. The latter then gives $\models^\alpha_{\mathfrak{C}^*, \ n} C$. Now, the induction hypothesis instantiated on $C$ yields 
            \begin{center} 
            $\forall \mathfrak{C} \supseteq_n \mathfrak{B} \ (\models^\alpha_{\mathfrak{C}, \ n} C \Longrightarrow \ \models^s_{\mathfrak{C}, \ n} C)$
            \end{center}
            which, when instantiated on $\mathfrak{C}^* \supseteq_n \mathfrak{C} \supseteq_n \mathfrak{B}$, yields $\models^s_{\mathfrak{C}^*, \ n} C$. So we obtain
            \begin{center} 
            $\forall \mathfrak{C}^* \supseteq_n \mathfrak{C} \ (\models^s_{\mathfrak{C}^*, \ n} B \Longrightarrow \ \models^s_{\mathfrak{C}^*, \ n} C)$
            \end{center}
            which, by point 2 of Definition 8, is equivalent to $B \models^s_{\mathfrak{C}, \ n} C$ which in turn, by point (d) of Definition 8, is equivalent to $\models^s_{\mathfrak{C}, \ n} B \rightarrow C$. Thus we have
            \begin{center} 
            $\forall \mathfrak{C} \supseteq_n \mathfrak{B} \ (\models^\alpha_{\mathfrak{C}, \ n} B \rightarrow C \Longrightarrow \ \models^s_{\mathfrak{C}, \ n} B \rightarrow C)$.
            \end{center}
    \end{itemize}
    \end{itemize}
    Suppose now $\Gamma \neq \emptyset$. Take any arbitrary $\mathfrak{C} \supseteq_n \mathfrak{B}$ and suppose $\Gamma \models^\alpha_{\mathfrak{C}, \ n} A$. By the admissibility clause, this means
    \begin{center} 
    $\forall \mathfrak{C}^* \supseteq_n \mathfrak{C} \ (\models^\alpha_{\mathfrak{C}^*, \ n} \Gamma \Longrightarrow \ \models^\alpha_{\mathfrak{C}^*, \ n} A)$.
    \end{center} 
    Take any arbitrary $\mathfrak{C}^* \supseteq_n \mathfrak{C}$ and assume $\models^s_{\mathfrak{C}^*, \ n} \Gamma$. Let (1) be instantiated on $\emptyset$, the elements of $\Gamma$, and $\mathfrak{C}^* \supseteq_n \mathfrak{C} \supseteq_n \mathfrak{B}$, so
    \begin{center} 
    $\models^s_{\mathfrak{C}^*, \ n} \Gamma \Longrightarrow \ \models^\alpha_{\mathfrak{C}^*, \ n} \Gamma$
    \end{center}
    which yields $\models^\alpha_{\mathfrak{C}^*, \ n} \Gamma$. The latter, when 
    \begin{center} 
    $\forall \mathfrak{C}^* \supseteq_n \mathfrak{C} \ (\models^\alpha_{\mathfrak{C}^*, \ n} \Gamma \Longrightarrow \ \models^\alpha_{\mathfrak{C}^*, \ n} A)$
    \end{center}
    is instantiated on $\mathfrak{C}^*$, yields $\models^\alpha_{\mathfrak{C}^*, \ n} A$. By what proved for the case with $\Gamma = \emptyset$, we obtain that
    \begin{center} 
    $\forall \mathfrak{C}^* \supseteq_n \mathfrak{C} \ (\models^s_{\mathfrak{C}^*, \ n} \Gamma \Longrightarrow \ \models^s_{\mathfrak{C}^*, \ n} A)$,
    \end{center}
    which means $\Gamma \models^s_{\mathfrak{C}, \ n} A$ by point 2 of Definition 8.
\end{proof}

\begin{corollary}
    $\forall \Gamma \ \forall A \ \forall \mathfrak{B} \in \mathbb{B}^n \ (\Gamma \models^s_{\mathfrak{B}, \ n} A \Longrightarrow \ \Gamma \models^\alpha_{\mathfrak{B}, \ n} A) \Longrightarrow \ \forall \Gamma \ \forall A \ \forall \mathfrak{B} \in \mathbb{B}^n \ (\Gamma \models^\alpha_{\mathfrak{B}, \ n} A \Longrightarrow \ \Gamma \models^s_{\mathfrak{B}, \ n} A)$.\footnote{If we want to stick to the ‘‘deviant" version of reducibility semantics of footnote 3 we have similar results.

\begin{theorem}
    $\forall \Gamma \ \forall A \ \forall \mathfrak{C} \supseteq_n \mathfrak{B} \ (\Gamma \models_{\mathfrak{C}, \ n} A \Longrightarrow \Gamma \models^\alpha_{\mathfrak{C}, \ n} A) \Longrightarrow \forall \Gamma \ \forall A \ \forall \mathfrak{C} \supseteq_n \mathfrak{B} \ (\Gamma \models^\alpha_{\mathfrak{C}, \ n} A \Longrightarrow \Gamma \models_{\mathfrak{C}, \ n} A)$.
\end{theorem}

\begin{corollary}
    $\forall \Gamma \ \forall A \ \forall \mathfrak{B} \in \mathbb{B}^n \ (\Gamma \models_{\mathfrak{B}, \ n} A \Longrightarrow \Gamma \models^\alpha_{\mathfrak{B}, \ n} A) \Longrightarrow \forall \Gamma \ \forall A \ \forall \mathfrak{B} \in \mathbb{B}^n \ (\Gamma \models^\alpha_{\mathfrak{B}, \ n} A \Longrightarrow \Gamma \models_{\mathfrak{B}, \ n} A).$
\end{corollary}}
\end{corollary}

\begin{proposition}
    On a disjunction-free propositional language, $\Gamma \models^\alpha_{\mathfrak{B}, \ n} A \Longleftrightarrow \Gamma \models^s_{\mathfrak{B}, \ n} A$.
\end{proposition}

\begin{proposition}
    On a disjunction-free propositional language, $\Gamma \models^\alpha_n A \Longleftrightarrow \Gamma \models^s_n A$.
\end{proposition}

\subsection{Base-incomparability}

The different ways in which Prawitz's reducibility semantics and Sandqvist's base semantics deal with $\vee$ do not exclude by themselves that the two approaches might be somehow ‘‘ordered" relative to consequence over an atomic base. Let us say that $\models^\alpha$ and $\models^s$ are \emph{base-comparable} if and only if either the antecedent or the consequent of Corollary 3 hold. Base-comparability can be also read as the property that $\models^\alpha$ and $\models^s$ are ‘‘models-monomorphic", i.e., either all models of $\models^\alpha$ are also models of $\models^s$, or vice versa, i.e., again, more precisely, if we set

\begin{center}
    $\mathbb{M}^\alpha_{\Gamma, \ A, \ n} = \{\mathfrak{B} \ | \ \Gamma \models^\alpha_{\mathfrak{B}, \ n} A\}$ and $\mathbb{M}^s_{\Gamma, \ A, \ n} = \{\mathfrak{B} \ | \ \Gamma \models^s_{\mathfrak{B}, \ n} A\}$,
\end{center}
then for all $n, \Gamma$ and $A$ either $\mathbb{M}^\alpha_{\Gamma, \ A, \ n} \subseteq \mathbb{M}^s_{\Gamma, \ A, \ n}$ or $\mathbb{M}^s_{\Gamma, \ A, \ n} \subseteq \mathbb{M}^\alpha_{\Gamma, \ A, \ n}$.

However, given what we know today at the level of \emph{logical validity}, we can also immediately rule out that, in some (very relevant) cases, Prawitz's reducibility semantics is ‘‘models-monomorphic" over Sandqvist's base semantics. For, $\texttt{IL}$ is known to be complete over the latter when atomic bases have level $\geq 2$---see \cite{sandqvist} and Theorem 9 below---while we have said above that $\texttt{IL}$ is \emph{never} complete over the former---see Theorem 3 above. Therefore, since logical validity means validity over all atomic bases, the consequent of Theorem 4 and Corollary 3 fails when atomic bases have level $\geq 2$. And now Theorem 4 and Corollary 3 come to play an active role by themselves, for via them we can infer in turn, by contraposition, that under the given conditions their antecedent fails too, namely, that Sandqvist's base semantics is not ‘‘models-monomorphic" over Prawitz's reducibility semantics either with atomic bases of level $\geq 2$. So, under the given conditions, Prawitz's reducibility semantics and Sandqvist's base semantics are not base-comparable or, to put it in another way, their classes of models $\mathbb{M}^\alpha_{\Gamma, \ A, \ n}$ and $\mathbb{M}^s_{\Gamma, \ A, \ n}$ diverge when $n \geq 2$. In turn, this may have connections with Schroeder-Heister's proposal to consider Sandqvist's approach as belonging to a family of proof-theoretic semantics which prioritise elimination rules \emph{in general} \cite{schroederheisterrolf}---for an elimination-based approach, see also \cite{GheorghiuPym, hermogenespragmatist}. This interesting topic is currently under investigation in PTS, but I shall leave it aside here.\footnote{This was pointed out to me by Peter Schroeder-Heister, whom I am therefore indebted to for what follows in this section. Oliveira's approach abstracts from the kind of atomic bases, while the interest of the comparison at issue here lies in the fact that the approaches use the \emph{same} bases.}

As we shall see, however, it is part of the second aforementioned consequence of Theorem 4 and Corollary 3 that the latter can be also used to get some \emph{positive} information. Before turning to that, let me first make precise the informal remarks that I have been carrying out thus far. For doing this, I must refer to some concepts and results presented in \cite{piechaschroeder-heisterdecampossanz, piechaschroeder-heister2019} and \cite{sandqvist}---some of which will be used also in Section 5. Piecha, de Campos Sanz and Schroeder-Heister \cite{piechaschroeder-heisterdecampossanz} observed that any disjunction-free formula $A$ can be associated to a set of atomic rules. The association requires first of all to transform such an $A$ into a suitable $\#(A)$ via a number of step-wise replacements---with respect to which $\#(A)$ is irreducible---as follows: any sub-formula of the form $B \rightarrow C \wedge D$ is replaced by $(B \rightarrow C) \wedge (B \rightarrow D)$, and any sub-formula of the form $B \rightarrow (C \rightarrow D)$ is replaced by $(B \wedge C) \rightarrow D$---note that, thereby, the consequent of every implication becomes atomic.

\begin{definition}[Piecha, de Campos Sanz \& Schroeder-Heister \cite{piechaschroeder-heisterdecampossanz}]
    To any disjunction-free $A$ we associate a set of atomic rules via a function $^\circ$ defined as follows:
    \begin{itemize}
        \item $\#(A) \in \texttt{ATOM}_\mathscr{L} \Longrightarrow A^\circ = \{A\}$;
        \item $\#(A) = B_1 \wedge ... \wedge B_n \rightarrow C$ with $C \in \texttt{ATOM}_\mathscr{L} \Longrightarrow A^\circ = \{R\}$, where $R$ has the form
        \begin{prooftree}
            \AxiomC{$[\Re_1]$}
            \noLine
            \UnaryInfC{$D_1$}
            \AxiomC{}
            \noLine
            \UnaryInfC{$\dots$}
            \AxiomC{$[\Re_n]$}
            \noLine
            \UnaryInfC{$D_n$}
            \TrinaryInfC{$C$}
        \end{prooftree}
        and
        \begin{prooftree}
            \AxiomC{$\Re_1$}
            \noLine
            \UnaryInfC{$D_1$}
            \AxiomC{}
            \noLine
            \UnaryInfC{$\dots$}
            \AxiomC{$\Re_n$}
            \noLine
            \UnaryInfC{$D_n$}
            \noLine
            \TrinaryInfC{}
        \end{prooftree}
        correspond to $B^\circ_1, ..., B^\circ_n$ respectively;
        \item $\#(A) = B_1 \wedge ... \wedge B_n \Longrightarrow A^\circ = \{B^\circ_1, ..., B^\circ_n\}$.
    \end{itemize}
\end{definition}

\begin{definition}
    Given a disjunction-free $\Gamma$, we set $\Gamma^\circ = \bigcup_{A \in \Gamma} A^\circ$.
\end{definition}

\noindent Let us now drop for a moment the constraint that the level of the atomic bases has an upper bound $n$. This notion can be easily obtained from Definitions 7, 8, 9, 20 and 21, by removing the requirement that the level of the atomic base is $m \leq n$, and that the extensions of the base have level at most $n$. I will indicate this notion by $\models$, while the class of atomic bases with no upper bound will be written $\mathbb{B}$. The extension-relation will be just $\supseteq$ (this is the only point where I use atomic bases with unlimited complexity, although I come back to this in the concluding remarks). Let me now state what I shall call \emph{import principle}---inspired by a principle called \emph{Import} which, in \cite{piechaschroeder-heister2019}, in turn based on \cite{piechaschroeder-heisterdecampossanz}, is stated for a number of frameworks and used to prove a number of results.

\begin{definition}
  $\models$ \emph{enjoys the import principle} iff, for every $\mathfrak{B} \in \mathbb{B}$ and every disjunction-free $\Gamma$, there is a disjunction-free $\Delta$ such that $\models_{\mathfrak{B} \cup \Delta^\circ} A \Longleftrightarrow \Gamma \models_\mathfrak{B} A$.
\end{definition}

\noindent Let us now consider what Piecha and Schroeder-Heister, in \cite{piechaschroeder-heister2019}, call the \emph{generalised disjunction property} (shortly, GDP): for every $\mathfrak{B} \in \mathbb{B}$, if $\Gamma$ is disjunction-free, then

\begin{center}
    $\Gamma \models_\mathfrak{B} A \vee B \Longrightarrow (\Gamma \models_\mathfrak{B} A$ or $\Gamma \models_\mathfrak{B} B)$.
\end{center}
I will now appeal to some notions and results concerning soundness and completeness of (recursive) systems---and of $\texttt{IL}$ in particular---over proof-theoretic validity. Both these notions and results will be more thoroughly explained in the next section. We know that, for every $n$, $\texttt{IL}$ is sound with respect to $\models_n$. This also holds in the case when the level of atomic bases has no upper bound. Using this fact, the following can be proved.

\begin{theorem}[Piecha \& Schroeder-Heister, \cite{piechaschroeder-heister2019}]
    If GDP holds on every $\mathfrak{B} \in \mathbb{B}$, Harrop's rule is valid in $\models$ (which, since Harrop's rule is not derivable in $\emph{\texttt{IL}}$ \cite{harrop}, means that the latter is then incomplete over $\models$).
\end{theorem}

\begin{proof}
    Take any arbitrary $\mathfrak{B} \in \mathbb{B}$, and $\mathfrak{C} \supseteq \mathfrak{B}$. Suppose $\models_\mathfrak{C} \neg A \rightarrow B \vee C$. Then, $\neg A \models_\mathfrak{C} B \vee C$. It is well-known that there is a disjunction-free $A^*$ such that $\neg A \vdash_{\texttt{IL}} A^*$ and $A^* \vdash_{\texttt{IL}} \neg A$. So, by soundness of $\texttt{IL}$, $A^* \models_\mathfrak{C} B \vee C$. By GDP, $A^* \models_\mathfrak{C} B$ or $A^* \models_\mathfrak{C} C$ and, again by soundness of $\texttt{IL}$, $\neg A \models_\mathfrak{C} B$ or $\neg A \models_\mathfrak{C} C$. Hence, $\models_\mathfrak{C} \neg A \rightarrow B$ or $\models_\mathfrak{C} \neg A \rightarrow C$, so $\models_\mathfrak{C} (\neg A \rightarrow B) \vee (\neg A \rightarrow C)$. Quantify over every $\mathfrak{C} \supseteq \mathfrak{B}$, and every $\mathfrak{B} \in \mathbb{B}$.
\end{proof}

\begin{theorem}[Piecha \& Schroeder-Heister \cite{piechaschroeder-heister2019}]
    If $\models$ enjoys the import principle, then GDP holds on every $\mathfrak{B} \in \mathbb{B}$.
\end{theorem}

\begin{proof}
    Assume $\Gamma \models_\mathfrak{B} A \vee B$ for $\Gamma$ disjunction-free. By import, there is a disjunction-free $\Delta$ such that $\models_{\mathfrak{B} \cup \Delta^\circ} A \vee B$, so $\models_{\mathfrak{B} \cup \Delta^\circ} A$ or $\models_{\mathfrak{B} \cup \Delta^\circ} B$ whence, again by import, $\Gamma \models_\mathfrak{B} A$ or $\Gamma \models_\mathfrak{B} B$.
\end{proof}

\begin{corollary}[Piecha \& Schroeder-Heister \cite{piechaschroeder-heister2019}]
    If $\models$ enjoys the import principle, then $\texttt{IL}$ is incomplete over $\models$.
\end{corollary}

Piecha, de Campos Sanz and Schroeder-Heister also proved that $\models$ does enjoy the import principle---by taking $\Delta$ to be just $\Gamma$, which gives $\Gamma^\circ$ as defined above---and that, via Sanqvist's coding \cite{sandqvist}, atomic rules of level $\geq 3$ can be reduced to atomic rules of level $2$ \cite{piechaschroeder-heisterdecampossanz}.

\begin{theorem}[Piecha, de Campos Sanz \& Schroeder-Heister \cite{piechaschroeder-heisterdecampossanz}]
    $\texttt{IL}$ is incomplete over $\models_2$.
\end{theorem}

\noindent Theorem 8 implies Theorem 3 above. On the other hand, we have what follows.

\begin{theorem}[Sandqvist \cite{sandqvist}]
$\texttt{IL}$ is complete over $\models^s_2$.\footnote{As said, Sandqvist's completeness proof applies to a language where $\bot$ is a nullary logical constant equipped with a special clause stating that $\bot$ holds on a base iff every atom holds on the base. The proof can be adapted to a language where $\bot$ is instead an atomic constant, and where bases come with atomic explosion. Via Proposition 2 above, what one obtains in this way is just Sandqvist's special clause for $\bot$. Then, in the construction of the ‘‘tailored" base for any valid sequent, $\bot$ is mapped onto itself. The ‘‘tailored" base will at that point contain, by default, a rule which infers any atomic image under such a mapping from derivations of (the image of) $\bot$ (under the mapping). For further details see \cite{sandqvist}.}
\end{theorem}

\begin{theorem}
    With $n = 2$, both the consequent and the antecedent of Corollary 3 fail.
\end{theorem}

\begin{proof}
    By Theorem 3, for some $\Gamma$ and $A$, $\Gamma \models^\alpha_2 A$ and $\Gamma \not\vdash_{\texttt{IL}} A$. Now $\Gamma \models^\alpha_2 A$ implies that $\forall \mathfrak{B} \in \mathbb{B}^2 \ (\Gamma \models_{\mathfrak{B}, \ 2} A)$. If the consequent of Corollary 3 holds on $n = 2$, we have $\forall \mathfrak{B} \in \mathbb{B}^2 \ (\Gamma \models^s_{\mathfrak{B}, \ 2} A)$, i.e. $\Gamma \models^s_2 A$. But Theorem 9 implies $\Gamma \not\models^s_2 A$. So, the consequent of Corollary 3 fails on $n = 2$. Hence, the antecedent of Corollary 3 fails on $n = 2$ too.
\end{proof}

\noindent If we allow for classical logic in the meta-language, we can now appeal to Theorem 10 to refine the informal description provided above, about the first relevant consequence of Theorem 4 and Corollary 3. I.e., we can infer from Theorem 10 (by classical meta-logic) that there are $\Gamma, A$ and $\mathfrak{B} \in \mathbb{B}^2$ such that $\Gamma \models^\alpha_{\mathfrak{B}, \ 2} A$ and $\Gamma \not\models^s_{\mathfrak{B}, \ 2} A$, and $\Gamma^*, A^*$ and $\mathfrak{B}^* \in \mathbb{B}^2$ such that $\Gamma^* \models^s_{\mathfrak{B}^*, \ 2} A^*$ and $\Gamma^* \not\models^\alpha_{\mathfrak{B}^*, \ 2} A^*$. As a final observation, let me also stress that, via Theorems 1 and 2, the comparative results above hold also for a comparison between standard base semantics and Sanqvist's reading.

\begin{theorem}
    $\forall \Gamma \ \forall A \ \forall \mathfrak{C} \supseteq_n \mathfrak{B} \ (\Gamma \models^s_{\mathfrak{C}, \ n} A \Longrightarrow \Gamma \models_{\mathfrak{C}, \ n} A) \Longrightarrow \forall \Gamma \ \forall A \ \forall \mathfrak{C} \supseteq_n \mathfrak{B} \ (\Gamma \models_{\mathfrak{C}, \ n} A \Longrightarrow \Gamma \models^s_{\mathfrak{C}, \ n} A)$.
\end{theorem}

\begin{corollary}
    $\forall \Gamma \ \forall A \ \forall \mathfrak{B} \in \mathbb{B}^n \ (\Gamma \models^s_{\mathfrak{B}, \ n} A \Longrightarrow \Gamma \models_{\mathfrak{B}, \ n} A) \Longrightarrow \forall \Gamma \ \forall A \ \forall \mathfrak{B} \in \mathbb{B}^n \ (\Gamma \models_{\mathfrak{B}, \ n} A \Longrightarrow \Gamma \models^s_{\mathfrak{B}, \ n} A)$.
\end{corollary}

\begin{theorem}
    With $n = 2$, both the antecedent and the consequent of Corollary 6 fail.\footnote{The same applies to the results mentioned above, concerning the disjunction-free propositional language.
    
    \begin{proposition}
    On a disjunction-free propositional language, $\Gamma \models_{\mathfrak{B}, \ n} A \Longleftrightarrow \Gamma \models^s_{\mathfrak{B}, \ n} A$.
\end{proposition}

\begin{proposition}
    On a disjunction-free propositional language, $\Gamma \models_n A \Longleftrightarrow \Gamma \models^s_n A$.
\end{proposition}
\noindent The proofs of Theorems 11 and 12, Corollaries 7, and Propositions 9 and 10 are the same as in the original case.}

\end{theorem}

\noindent A concrete example of the failure of the antecedents of Theorems 10 and 12 is the atomic base $\{R\}$ of Section 2, where $R$ is the (schematic) rule

\begin{prooftree}
    \AxiomC{$A$}
    \AxiomC{$[B]$}
    \noLine
    \UnaryInfC{$D$}
    \AxiomC{$[C]$}
    \noLine
    \UnaryInfC{$D$}
    \TrinaryInfC{$D$}
\end{prooftree}
for some $A, B, C \in \texttt{ATOM}_\mathscr{L}$ and every $D \in \texttt{ATOM}_\mathscr{L}$. We have $A \models^s_{\{R\}, \ 2} B \vee C$, but $A \not\models^\alpha_{\{R\}, \ 2} B \vee C$ and $A \not\models_{\{R\}, \ 2} B \vee C$.

\section{On base-completeness}

Besides implying that Prawitz's reducibility semantics and Sandqvist's base semantics are not base-comparable, Theorem 4 and Corollary 3 have as said also a second, ‘‘positive" implication. They provide a sufficient condition for both Prawitz's reducibility semantics and Sandqvist's base semantics to be equivalent relative to \emph{logical} consequence, and for a logic to be complete over Prawitz's reducibility semantics. To show this, I introduce notions of ‘‘point-wise" soundness and completeness of given logics over proof-theoretic semantics. The interest of these notions, however, does not stem from interactions with Theorem 4 and Corollary 3 only, since the very same (in)completeness phenomena that motivated the negative consequences of these results, also imply limit-results as concerns ‘‘point-wise" completeness of given logics over proof-theoretic semantics in general---in particular, for $\texttt{IL}$.

\subsection{Basic definitions and results}

By a \emph{system} $\Sigma$ I shall understand a recursive set of super-intuitionistic rules (over $\mathscr{L}$). The derivability of $A$ from $\Gamma$ in $\Sigma$ is indicated as usual with the notation $\Gamma \vdash_\Sigma A$. For example, $\Sigma$ may be $\texttt{IL}$. To deal with derivability at the atomic level, $\Sigma$ may be required to incorporate rules from an atomic base $\mathfrak{B}$. This yields the following general definition.

\begin{definition}
    The \emph{extended derivations-set} $\texttt{DER}^+_{\Sigma}$ of $\Sigma$ is defined inductively as follows:

    \begin{itemize}
        \item the single node labelled by $A \in \texttt{FORM}_\mathscr{L}$, possibly used as an axiom when $A \in \texttt{ATOM}_\mathscr{L}$, is a derivation in $\texttt{DER}^+_{\Sigma}$;
        \item if the following is a derivation in $\texttt{DER}^+_{\Sigma}$,
        \begin{prooftree}
            \AxiomC{$\Gamma_i, \mathfrak{C}_i, \Re_i$}
            \noLine
            \UnaryInfC{$\mathscr{D}_i$}
            \noLine
            \UnaryInfC{$A_i$}
        \end{prooftree}
        where $\mathfrak{C}_i, \Re_i$ are sets of atomic rules used in $\mathscr{D}_i$ and $A_i$ is the premise of an atomic rule $R$ of the form
        \begin{prooftree}
            \AxiomC{$[\Re_1]$}
            \noLine
            \UnaryInfC{$A_1$}
            \AxiomC{$\dots$}
            \AxiomC{$[\Re_n]$}
            \noLine
            \UnaryInfC{$A_n$}
            \RightLabel{$R$}
            \TrinaryInfC{$B$}
        \end{prooftree}
        ($i \leq n$)---where it is not necessarily required that $R \in \mathfrak{B}$---then
        \begin{prooftree}
            \AxiomC{$\Gamma_1, \mathfrak{C}_1, [\Re_1]$}
            \noLine
            \UnaryInfC{$\mathscr{D}_1$}
            \noLine
            \UnaryInfC{$A_1$}
            \AxiomC{$\dots$}
            \AxiomC{$\Gamma_n, \mathfrak{C}_n, [\Re_n]$}
            \noLine
            \UnaryInfC{$\mathscr{D}_n$}
            \noLine
            \UnaryInfC{$A_n$}
            \RightLabel{$R$}
            \TrinaryInfC{$B$}
        \end{prooftree}
        is a derivation in $\texttt{DER}^+_{\Sigma}$;
        \item the case of the logical rules of $\Sigma$ runs in a standard inductive way.
    \end{itemize}
\end{definition}
\noindent So, for example, when $\Sigma$ is $\texttt{IL}$, the last clause means that derivations in $\texttt{DER}^+_{\texttt{IL}}$ are defined by smooth induction when standard introduction and elimination rules are at issue.

\begin{definition}
    Let $\Re$ be a finite set of atomic rules such that $\Re \cap \mathfrak{B} = \emptyset$. That $A$ is \emph{derivable from} $\Gamma, \Re$ \emph{in} $\Sigma \cup \mathfrak{B}$ is indicated by $\Gamma, \Re \vdash_{\Sigma \cup \mathfrak{B}} A$, and it holds iff there is a derivation in $\texttt{DER}^+_{\Sigma}$ from $\Gamma$ to $A$ whose only additional rules besides those in $\mathfrak{B}$ are those in $\Re$.
\end{definition}
\noindent In the following, $\Vdash_n$ can indifferently be either $\models_n$, or $\models^\alpha_n$, or $\models^s_n$.

\begin{definition}
    $\Sigma$ is \emph{base-complete} over $\Vdash_n$ iff $\forall \mathfrak{B} \in \mathbb{B}^n, \Gamma \Vdash_{\mathfrak{B}, \ n} A \Longrightarrow \Gamma \vdash_{\Sigma \cup \mathfrak{B}} A$.
\end{definition}

\begin{definition}
    $\Sigma$ is \emph{base-sound} over $\Vdash_n$ iff $\forall \mathfrak{B} \in \mathbb{B}^n, \Gamma \vdash_{\Sigma \cup \mathfrak{B}} A \Longrightarrow \Gamma \Vdash_{\mathfrak{B}, \ n} A$.
\end{definition}

\begin{definition}
    $\Sigma$ is \emph{sound} over $\Vdash_n$  iff $\Gamma \vdash_\Sigma A \Longrightarrow \Gamma \Vdash_n A$.
\end{definition}

\begin{definition}
    $\Sigma$ is \emph{complete} over $\Vdash_n$ iff $\Gamma \Vdash_n A \Longrightarrow \Gamma \vdash_\Sigma A$.
\end{definition}

\begin{proposition}
    $\Sigma$ base-sound over $\Vdash_n \ \Longrightarrow \ \Sigma$ sound over $\Vdash_n$.
\end{proposition}

\begin{proof}
    Take $\Gamma \vdash_\Sigma A$. This means $\Gamma \vdash_{\Sigma \cup \mathfrak{B}^\emptyset} A$. By assumption of base-soundness of $\Sigma$, $\Gamma \Vdash_{\mathfrak{B}^\emptyset, \ n} A$, i.e. $\Gamma \Vdash_n A$.
\end{proof}

\begin{proposition}
    $\Sigma$ base-complete over $\Vdash_n \ \Longrightarrow \ \Sigma$ complete over $\Vdash_n$.
\end{proposition}

\begin{proof}
    Take $\Gamma \Vdash_n A$. This means $\Gamma \Vdash_{\mathfrak{B}, \ n} A$ for every $\mathfrak{B} \in \mathbb{B}^n$ and, by assumption of base-completeness of $\Sigma$, for every $\mathfrak{B} \in \mathbb{B}^n$, $\Gamma \vdash_{\Sigma \cup \mathfrak{B}} A$. By instantiating this on $\mathfrak{B}^\emptyset$, $\Gamma \vdash_{\Sigma \cup \mathfrak{B}^\emptyset} A$, i.e., $\Gamma \vdash_\Sigma A$.
\end{proof}

\noindent From this, we can immediately infer base-incompleteness of $\texttt{IL}$ over $\Vdash_n$, where $\Vdash_n$ is now either $\models_n$ or $\models^\alpha_n$.

\begin{proposition}
    $\texttt{IL}$ is not base-complete over $\Vdash_n$.
\end{proposition}

\subsection{Applications to the comparison}

\noindent We can now use the notions of base-soundness and base-completeness to single out a sufficient condition for reducibility semantics and Sandqvist's base semantics to be equivalent.

\begin{theorem}
    $\exists \Sigma \ (\Sigma$ base-complete over $\models^s_n$ and base-sound over $\models^\alpha_n) \Longrightarrow \forall \Gamma \ \forall A \ (\Gamma \models^s_n A \Longleftrightarrow \Gamma \models^\alpha_n A)$.
\end{theorem}

\begin{proof}
   Let $\Sigma$ be as required, and take an arbitrary $\Gamma$ and $A$. Let us prove the direction ($\Longrightarrow$). Assume $\Gamma \models^s_n A$. By assumption of base-completeness, $\Gamma \vdash_\Sigma A$. Again, by assumption of base-soundness, $\Gamma \models^\alpha_n A$. Let us prove the direction ($\Longleftarrow$). Take arbitrary $\Delta$, $B$ and $\mathfrak{B} \in \mathbb{B}^n$, and assume $\Delta \models^s_{\mathfrak{B}, \ n} B$. By assumption of base-completeness, we have $\Delta \vdash_{\Sigma \cup \mathfrak{B}} B$ and, by assumption of base-soundness, we have $\Delta \models^\alpha_{\mathfrak{B}, \ n} B$. Hence,
   \begin{center}
   $\Delta \models^s_{\mathfrak{B}, \ n} B \Longrightarrow \Delta \models^\alpha_{\mathfrak{B}, \ n} B$.
   \end{center}
   We can now introduce universal quantification over $\Delta, B$ and $\mathfrak{B} \in \mathbb{B}^n$, so to obtain
   \begin{center}
       $\forall \Delta \ \forall B \ \forall \mathfrak{B} \in \mathbb{B}^n \ (\Delta \models^s_{\mathfrak{B}, \ n} B \Longrightarrow \Delta \models^\alpha_{\mathfrak{B}, \ n} B)$.
   \end{center}
   By Corollary 3, this implies
   \begin{center}
       $\forall \Delta \ \forall B \ \forall \mathfrak{B} \in \mathbb{B}^n \ (\Delta \models^\alpha_{\mathfrak{B}, \ n} B \Longrightarrow \Delta \models^s_{\mathfrak{B}, \ n} B)$.
   \end{center}
   By instantiating on our previously chosen arbitrary $\Gamma$ and $A$, we obtain
   \begin{center} 
   $\forall \mathfrak{B} \in \mathbb{B}^n \ (\Gamma \models^\alpha_{\mathfrak{B}, \ n} A \Longrightarrow \Gamma \models^s_{\mathfrak{B}, \ n} A)$.
   \end{center}
   Assume now $\Gamma \models^\alpha_n A$. This means $\forall \mathfrak{B} \in \mathbb{B}^n \ (\Gamma \models^\alpha_{\mathfrak{B}, \ n} A)$, whence $\forall \mathfrak{B} \in \mathbb{B}^n \ (\Gamma \models^s_{\mathfrak{B}, \ n} A)$, which in turn means $\Gamma \models^s_n A$.
\end{proof}

\begin{corollary}
    $\Sigma$ base-complete over $\models^s_n$ and base-sound over $\models^\alpha_n \ \Longrightarrow \ \Sigma$ complete over $\models^\alpha_n$.
\end{corollary}

\begin{proof}
    Take any arbitrary $\Sigma$ as required and suppose $\Gamma \models^\alpha_n A$. Theorem 13 yields $\Gamma \models^s_n A$ which, by base-completeness of $\Sigma$ over $\models^s_n$, implies $\Gamma \vdash_\Sigma A$.\footnote{Via Theorems 1 and 2, Theorem 13 and Corollary 7 can be also formulated for standard base semantics.
    
\begin{theorem}
        $\exists \Sigma \ (\Sigma$ base-complete over $\models^s_n$ and base-sound over $\models_n) \Longrightarrow \forall \Gamma \ \forall A \ (\Gamma \models^s_n A \Longleftrightarrow \Gamma \models_n A)$.
\end{theorem}
    
    \begin{corollary}
        $\Sigma$ base-complete over $\models^s_n$ and base-sound over $\models_n \Longrightarrow \ \Sigma$ complete over $\models_n$.
    \end{corollary}
    \noindent Let me also stress that in the ‘‘deviant" reading of reducibility semantics mentioned in footnote 3---i.e., with only the left-to-right direction in the admissibility clause---something similar to Theorem 13 and Corollary 7 can be obtained for a relation between $\models_n$ and $\models^\alpha_n$.
    
    \begin{theorem}
        $\exists \Sigma \ (\Sigma$ base-complete over $\models_n$ and base-sound over $\models^\alpha_n) \Longrightarrow \forall \Gamma \ \forall A \ (\Gamma \models_n A \Longleftrightarrow \Gamma \models^\alpha_n A)$.
    \end{theorem}
    
    \begin{corollary}
        $\Sigma$ base-complete over $\models_n$ and base-sound over $\models^\alpha_n \ \Longrightarrow \ \Sigma$ complete over $\models^\alpha_n$.
    \end{corollary}
    \noindent The proofs of Theorem 15 and Corollary 9 are similar to those of Theorem 13 and Corollary 7, using the results mentioned in footnote 5.}
\end{proof}
\noindent Besides their conceptual interest (and their link with the results proved below), Theorem 13 and Corollary 7 might have a number of interesting applications. In the next section, I discuss more extensively the notion of base-completeness, which is involved in both the results at issue.

\subsection{Inconsistency of intuitionistic base-completeness (and beyond)}

Contrarily to model-theory, where the models of one's language are given by mappings from the language onto (mostly set-theoretic) structures, the ‘‘models" of proof-theoretic semantics are deductive in nature, that is, they are just proof-systems. In a sense---above all if one looks at Prawitz's first semantic works \cite{prawitz1971, prawitz1973}---proof-theoretic semantics might be even understood as a sort of ‘‘semantic generalisation" of normalisation theory for Natural Deduction, namely, as a semantics where certain normalisation properties \emph{provable} relative to given Natural Deduction systems, are turned into semantic \emph{requirements} which determine the meaning of the logical terminology, and which provide validity criteria for argument-structures. An example of this is what Schroeder-Heister called the \emph{fundamental corollary} of normalisation theory \cite{schroeder-heister2006}, stating that closed normal derivations in constructive systems end in introduction form. In proof-theoretic semantics, this becomes the tenet that a closed argument structure is valid (on an atomic base) when it reduces (relative to some set of reductions) to a canonical form whose immediate substructures are also valid (on the same atomic base and relative to the same set of reductions)---for more on this, see \cite{ptssquare, piccolominiattitudes}. More in general, this complies with the idea that meaning and validity are essentially embedded into, or even stemming from, deduction and structural features of deduction.

From this point of view, it may be natural to expect that logics which are sound or complete over a given variant $\Vdash$ of proof-theoretic consequence, are also base-sound or base-complete over $\Vdash$. Given the intertwinement between meaning and validity, on the one hand, and deduction and its structural properties, on the other hand, it may be in other words natural to expect that, if a logic $\Sigma$ is sound over $\Vdash$, derivations in $\Sigma$ from $\Gamma$ to $A$ can be turned into argument-structures from $\Gamma$ to $A$ which are valid relative to $\Vdash_n$ (and which are extracted from $\Gamma \Vdash_n A$, when $\Vdash$ is $\models$ or $\models^s$). Since derivations are invariant under addition of atomic bases $\mathfrak{B} \in \mathbb{B}^n$ to $\Sigma$, one may expect the same to happen when the corresponding valid argument-structures are evaluated over $\Vdash_{\mathfrak{B}, \ n}$. Conversely, if $\Sigma$ is complete over $\Vdash$, one may expect that argument-structures from $\Gamma$ to $A$ which are valid relative to $\Vdash$ (extracted from $\Gamma \Vdash_n A$, when $\Vdash$ is $\models$ or $\models^s$), can be ‘‘represented" by derivations in $\Sigma$, and that this ‘‘representation" property is stable relative to validity over $\mathfrak{B} \in \mathbb{B}^n$, namely, that an argument-structure valid relative to $\Vdash_{\mathfrak{B}, \ n}$ has a ‘‘representative" derivation in $\Sigma \cup \mathfrak{B}$.

In what follows, I shall not concentrate on base-soundness and base-completeness in general, but on base-soundness and base-completeness of $\texttt{IL}$. The main aim is showing that, although \emph{prima facie} reasonable, base-completeness is an inconsistent notion when referred to $\texttt{IL}$---and for proving this, I shall need to use (base-)soundness of $\texttt{IL}$ too. Let us first of all state the following---where $\Vdash$ means again as before $\models_n$, or $\models^\alpha_n$, or $\models^s_n$.

\begin{proposition}
    $\texttt{IL}$ is base-sound over $\Vdash_n$.
\end{proposition}

\begin{proof}
    By induction on the length of the derivation of $\Gamma \vdash_{\texttt{IL} \cup \mathfrak{B}} A$. One reasons like in \cite[Theorem 4.1]{sandqvist} for elimination of disjunction on $\models^s_{\mathfrak{B}, \ n}$, and uses \cite[Theorem 3.1]{sandqvist} or \cite[Remark 3.5]{piechaschroeder-heister2019} for preservation of validity by atomic rules.
\end{proof}

\begin{corollary}
    $\texttt{IL}$ is sound over $\Vdash_n$.
\end{corollary}
\noindent Let me now introduce what Piecha, de Campos Sanz and Schroeder-Heister have called the \emph{export principle} \cite{piechaschroeder-heisterdecampossanz}---see also \cite{gheorghiupymnegationasfailure}. Before doing this, it is useful to provide the preliminary definition of translation-function from an atomic base $\mathfrak{B}$ into a set of disjunction-free formulas $\mathfrak{B}^*$.

\begin{definition}[Piecha, de Campos Sanz \& Schroeder-Heister, \cite{piechaschroeder-heisterdecampossanz}]
    To any rule $R$ of level $n$ we associate a disjunction-free formula via a function $^*$ defined by induction as follows:
    \begin{itemize}
        \item $\mathfrak{L}(R) = 0 \Longrightarrow R = A$ with $A \in \texttt{ATOM}_\mathscr{L}$. Then $R^* = A$;
        \item $\mathfrak{L}(R) = k + 1 \Longrightarrow R$ has the form
        \begin{prooftree}
            \AxiomC{$[\Re_1]$}
            \noLine
            \UnaryInfC{$A_1$}
            \AxiomC{}
            \noLine
            \UnaryInfC{$\dots$}
            \AxiomC{$[\Re_n]$}
            \noLine
            \UnaryInfC{$A_n$}
            \TrinaryInfC{$A$}
        \end{prooftree}
        where
        \begin{prooftree}
            \AxiomC{$\Re_1$}
            \noLine
            \UnaryInfC{$A_1$}
            \AxiomC{}
            \noLine
            \UnaryInfC{$\dots$}
            \AxiomC{$\Re_n$}
            \noLine
            \UnaryInfC{$A_n$}
            \noLine
            \TrinaryInfC{}
        \end{prooftree}
        are rules $R_1, ..., R_n$ of level $\leq k$. Then $R^* = \bigwedge_{i \leq n} R^*_i \rightarrow A$.
    \end{itemize}
\end{definition}

\begin{definition}[Piecha, de Campos Sanz \& Schroeder-Heister, \cite{piechaschroeder-heisterdecampossanz}]
    Given $\mathfrak{B} = \{R_1, ..., R_n\}$, we set $\mathfrak{B}^* = \{R^*_1, ..., R^*_n\}$.
\end{definition}

\noindent Piecha, de Campos Sanz and Schroeder-Heister's export principle was defined for standard consequence over atomic bases with unlimited complexity. However, their definition can be easily adapted, both to consequence over bases with upper-bounded level, and to consequence in Sandqvist's sense.

More importantly, it is now time I deal with an issue that I could postpone so far, but which becomes crucial when discussing the export principle relative to what is at stake here. As stated in Section 2, I have been assuming throughout this paper

\begin{itemize}
    \item that the $\Gamma$-s in the consequence statements are always finite, and
    \item that atomic bases can be infinite, in the sense of containing infinitely many atomic rules.
\end{itemize}
The first assumption was officially motivated by issues of compact monotonicity in Prawitz's version of proof-theoretic semantics pointed out by \cite{stafford1}. But another motivation is that the assumptions are \emph{both} required if one wants Sandqvist's completeness proof to go through---which is of course something I wanted, given that Sandqvist's completeness result is used in the proof of such results as Theorem 10.

Now, the export principle as provided by Piecha, de Campos Sanz and Schroeder-Heister (restricted to atomic bases with limited complexity) reads officially: for every $\mathfrak{B} \in \mathbb{B}^n$ there is disjunction-free $\Delta$ such that $\Gamma \Vdash_{\mathfrak{B}, \ n} A \Longleftrightarrow \Gamma, \Delta \Vdash_n A$. Below, however, we shall need to use this by requiring $\Delta$ to explicitly contain only elements of $\mathfrak{B}^*$, as the latter has been given in Definition 32. However, a formulation of the export principle of the form

\begin{center}
    for every $\mathfrak{B} \in \mathbb{B}^n, \Gamma \Vdash_{\mathfrak{B}, \ n} A \Longleftrightarrow \Gamma, \mathfrak{B}^* \Vdash_n A$
\end{center}
would not work. Since we allowed atomic bases to be infinite, $\mathfrak{B}^*$ might well be infinite too, but then the requirement on sets of assumptions in consequence statements to be finite would be violated. To get rid of this, I shall here adopt a stronger form of the export principle which, for reasons that will become clear in a while, I call the \emph{compact export principle}. This will be sufficient to prove a number of results. In what follows, $\Vdash_n$ stands again for $\models_n, \models^\alpha_n$ or $\models^s_n$.

\begin{definition}
    $\Vdash_n$ \emph{enjoys the compact export principle} iff $\forall \mathfrak{B} \in \mathbb{B}^n \ (\Gamma \Vdash_{\mathfrak{B}, \ n} A \Longleftrightarrow \exists \Delta \subseteq \mathfrak{B}^*$ finite such that $\Gamma, \Delta \Vdash_n A)$.
\end{definition}

\noindent The generalised disjunction property, used in Section 4 as referred to bases with unlimited level, can be now restricted to a consequence over bases with upper bound (shortly, $\text{GDP}^n$): for every $\mathfrak{B} \in \mathbb{B}^n$, if $\Gamma$ is disjunction-free, then

\begin{center}
    $\Gamma \Vdash_{\mathfrak{B}, \ n} A \vee B \Longrightarrow \ (\Gamma \Vdash_{\mathfrak{B}, \ n} A$ or $\Gamma \Vdash_{\mathfrak{B}, \ n} B)$.
\end{center}

\noindent A result similar to Theorem 6 can be proved for $\text{GDP}^n$. The proof runs in very much the same way as that of Theorem 6 itself, except that, when $\Vdash_n$ is $\models^s_n$, one has to use the following fact.

\begin{proposition}[Sandqvist \cite{sandqvist}]
    $(\models^s_{\mathfrak{B}, \ n} A$ or $\models^s_{\mathfrak{B}, \ n} B) \Longrightarrow \ \models^s_{\mathfrak{B}, \ n} A \vee B$.
\end{proposition}

\begin{theorem}
    If $\text{GDP}^n$ holds on every $\mathfrak{B} \in \mathbb{B}^n$, Harrop's rule is valid in $\Vdash_n$ so, since Harrop's rule is not derivable in $\emph{\texttt{IL}}$, the latter is incomplete over $\Vdash_n$.
\end{theorem}
   
\noindent Of course, since $\texttt{IL}$ is complete over $\models^s_2$, $GDP^2$ fails for $\models^s_2$. A counterexample is again the base $\{R\}$ mentioned in Section 2. We have $A \models^s_{\{R\}} B \vee C$, but neither $A \models^s_{\{R\}} B$, nor $A \models^s_{\{R\}} C$ (I am indebted to Hermógenes Oliveira for useful discussions on this topic). Next, we prove the following results---drawn from \cite{piechaschroeder-heister2019}.

\begin{theorem}
    Compact export principle plus completeness of $\texttt{IL}$ over $\Vdash_n$ imply that $\text{GDP}^n$ holds for every $\mathfrak{B} \in \mathbb{B}^n$.
\end{theorem}

\begin{proof}
    Suppose $\vee$ does not occur in $\Gamma$ and, for arbitrary $\mathfrak{B} \in \mathbb{B}^n$, assume $\Gamma \Vdash_{\mathfrak{B}, \ n} A \vee B$. Compact export implies that $\exists \Delta \subseteq \mathfrak{B}^*$ finite such that $\Gamma, \Delta \Vdash_n A \vee B$ and, by completeness of $\texttt{IL}$ over $\Vdash_n$, we have that $\exists \Delta \subseteq \mathfrak{B}^*$ finite such that $\Gamma, \Delta \vdash_{\texttt{IL}} A \vee B$. But $\vee$ does not occur in $\Delta$, so we have that $\exists \Delta \subseteq \mathfrak{B}^*$ finite such that $(\Gamma, \Delta \vdash_{\texttt{IL}} A$ or $\Gamma, \Delta \vdash_{\texttt{IL}} B)$, and this implies that $\exists \Delta \subseteq \mathfrak{B}^*$ finite such that $\Gamma, \Delta \vdash_{\texttt{IL}} A$, or $\exists \Delta \subseteq \mathfrak{B}^*$ finite such that $\Gamma, \Delta \vdash_{\texttt{IL}} B$. By soundness of $\texttt{IL}$, we have that $\exists \Delta \subseteq \mathfrak{B}^*$ finite such that $\Gamma, \Delta \Vdash_n A$, or $\exists \Delta \subseteq \mathfrak{B}^*$ finite such that $\Gamma, \Delta \Vdash_n B$. Again by compact export, $\Gamma \Vdash_{\mathfrak{B}, \ n} A$ or $\Gamma \Vdash_{\mathfrak{B}, \ n} B$.
\end{proof}

\begin{corollary}
    If $\Vdash_n$ enjoys the compact export principle, then $\texttt{IL}$ is incomplete over $\Vdash_n$.
\end{corollary}

\begin{proof}
By Theorem 17, compact export plus completeness of $\texttt{IL}$ over $\Vdash_n$ imply that $GDP^n$ holds for every $\mathfrak{B} \in \mathbb{B}^n$. By Theorem 16, the latter implies that Harrop’s rule is valid over $\Vdash_n$, hence that $\texttt{IL}$ is incomplete over $\Vdash_n$. Therefore, compact export plus completeness of $\texttt{IL}$ over $\Vdash_n$ imply incompleteness of $\texttt{IL}$ over $\Vdash_n$, that is, we can reject completeness of $\texttt{IL}$ over $\Vdash_n$ under the assumption that $\Vdash_n$ enjoys the compact export principle.
\end{proof}

\noindent Observe, in passing, that the latter implies what follows. 

\begin{proposition}
    $\models^s_2$ does not enjoy the compact export principle.
\end{proposition}

\noindent A counterexample is again given by the atomic base $\{R\}$ above. For, we have as said $A \models^s_{\{R\}} B \vee C$, but there can be no finite $\Delta \subset \{R\}^*$---i.e., no finite number of instances of $R$---such that $A, \Delta \models^s B \vee C$. Of course, this stems from the fact that, as per Sandqvist's disjunction clause, we need to quantify over every $D$ as occurring in $R$ for $B \vee C$ to hold.

The more general point, however, is that base-completeness of $\texttt{IL}$ at a given level is equivalent to compact export principle plus completeness of $\texttt{IL}$ at that level. Let us first establish a compact export principle for $\texttt{IL}$---see \cite[Remark 3.8]{piechaschroeder-heister2019} and \cite{schroederheisternaturalext} for a similar point. Starting with an example, let $\mathfrak{B}$ consist of the rules---leaving atomic explosion aside---

\begin{prooftree}
    \AxiomC{}
    \UnaryInfC{$p$}
    \AxiomC{$p$}
    \UnaryInfC{$v$}
    \AxiomC{$q$}
    \AxiomC{$r$}
    \BinaryInfC{$z$}
    \AxiomC{$[s]$}
    \noLine
    \UnaryInfC{$u$}
    \AxiomC{$v$}
    \BinaryInfC{$q$}
    \noLine
    \QuaternaryInfC{}
\end{prooftree}
and let $\Re$ consist of the rules

\begin{prooftree}
    \AxiomC{$s$}
    \RightLabel{$R_1$}
    \UnaryInfC{$t$}
    \AxiomC{}
    \RightLabel{$R_2$}
    \UnaryInfC{$r$}
    \noLine
    \BinaryInfC{}
\end{prooftree}
then the following

\begin{prooftree}
    \AxiomC{$[q \vee (t \rightarrow u)]^1$}
    \AxiomC{$[q]^2$}
    \AxiomC{}
    \LeftLabel{$R_2$}
    \UnaryInfC{$r$}
    \BinaryInfC{$z$}
    \AxiomC{$[s]^3$}
    \LeftLabel{$R_1$}
    \UnaryInfC{$t$}
    \AxiomC{$[t \rightarrow u]^4$}
    \BinaryInfC{$u$}
    \AxiomC{}
    \UnaryInfC{$p$}
    \UnaryInfC{$v$}
    \RightLabel{$3$}
    \BinaryInfC{$q$}
    \AxiomC{}
    \LeftLabel{$R_2$}
    \UnaryInfC{$r$}
    \BinaryInfC{$z$}
    \RightLabel{$2, 4$}
    \TrinaryInfC{$z$}
    \RightLabel{$1$}
    \UnaryInfC{$(q \vee (t \rightarrow u)) \rightarrow z$}
    \AxiomC{$w \rightarrow \bot$}
    \BinaryInfC{$((q \vee (t \rightarrow u)) \rightarrow z) \wedge (w \rightarrow \bot)$}
\end{prooftree}
witnesses that $w \rightarrow \bot, \Re \vdash_{\texttt{IL} \cup \mathfrak{B}} ((q \vee (t \rightarrow u)) \rightarrow z) \wedge (w \rightarrow \bot)$. Observe that this can be turned into

\begin{center}
$p, p \rightarrow v, w \rightarrow \bot, r, q \wedge r \rightarrow z, s \rightarrow t, ((s \rightarrow u) \wedge v) \rightarrow q \vdash_{\texttt{IL}} ((q \vee (t \rightarrow u)) \rightarrow z) \wedge (w \rightarrow \bot)$
\end{center}
where each assumption but $w \rightarrow \bot$ is $\rho^*$ with $\rho \in \Re \cup \mathfrak{B}$. This is what the following proposition establishes in general.

\begin{proposition}
    Let $\Re$ be a finite set of atomic rules. Then, $\Gamma, \Re \vdash_{\emph{\texttt{IL}} \cup \mathfrak{B}} A \Longleftrightarrow (\Gamma, \Re^*, \Delta \vdash_{\emph{\texttt{IL}}} A$, with $\Delta \subseteq \mathfrak{B}^*$ finite$)$.
\end{proposition}

\begin{proof}
    In order not to burden the notation excessively, I assume $\mathfrak{B}$ to be finite---hence I assume $\Delta = \mathfrak{B}^*$. The idea for the general case is, roughly, that only a \emph{finite} number of rules from $\mathfrak{B}$ is used in derivations in $\texttt{IL} \cup \mathfrak{B}$, thus the proof of this general case obtains from the one below via slight adjustments. ($\Longrightarrow$) We proceed by induction on the length $\lambda(\mathscr{D})$ of the derivation $\mathscr{D}$ in $\texttt{IL} \cup \mathfrak{B}$:
    \begin{itemize}
        \item $\lambda(\mathscr{D}) = 0 \Longrightarrow$ there are two cases. $\mathscr{D}$ is an application of the assumption rule, hence $A \in \Gamma$, hence $\Gamma \vdash_{\texttt{IL}} A$, hence $\Gamma, \Re^*, \mathfrak{B}^* \vdash_{\texttt{IL}} A$. Or $\mathscr{D}$ is an instance of an atomic rule $R$ of level $0$ in $\mathfrak{B}$ or in $\Re$, hence $A \in \texttt{ATOM}_\mathscr{L}$ and $A = R^* \in \mathfrak{B}^*$ or $A = R^* \in \Re$. But then clearly $\mathfrak{B}^* \vdash_{\texttt{IL}} A$, resp. $\Re^* \vdash_{\texttt{IL}} A$, hence $\Gamma, \Re^*, \mathfrak{B}^* \vdash_{\texttt{IL}} A$;
        \item $\lambda(\mathscr{D}) = k + 1 \Longrightarrow \mathscr{D}$ ends by
        \begin{itemize}
            \item application of ($\vee_E$) $\Longrightarrow \mathscr{D}$ has the form
            \begin{prooftree}
                \AxiomC{$\Gamma, \Re$}
                \noLine
                \UnaryInfC{$\mathscr{D}_1$}
                \noLine
                \UnaryInfC{$B \vee C$}
                \AxiomC{$[B], \Gamma, \Re$}
                \noLine
                \UnaryInfC{$\mathscr{D}_2$}
                \noLine
                \UnaryInfC{$A$}
                \AxiomC{$[C], \Gamma, \Re$}
                \noLine
                \UnaryInfC{$\mathscr{D}_3$}
                \noLine
                \UnaryInfC{$A$}
                \TrinaryInfC{$A$}
            \end{prooftree}
            So, $\Gamma, \Re \vdash_{\texttt{IL} \cup \mathfrak{B}} B \vee C$ via $\mathscr{D}_1$ with $\lambda(\mathscr{D}_1) \leq k$, $B, \Gamma, \Re \vdash_{\texttt{IL} \cup \mathfrak{B}} A$ via $\mathscr{D}_2$ with $\lambda(\mathscr{D}_2) \leq k$, and $C, \Gamma, \Re \vdash_{\texttt{IL} \cup \mathfrak{B}} A$ via $\mathscr{D}_3$ with $\lambda(\mathscr{D}_3) \leq k$. By induction hypothesis, there are $\mathscr{D}_1^*, \mathscr{D}_2^*$ and $\mathscr{D}_3^*$ in $\texttt{IL}$ such that $\Gamma, \Re^*, \mathfrak{B}^* \vdash_{\texttt{IL}} B \vee C$ holds via $\mathscr{D}_1^*$, $B, \Gamma, \Re^*, \mathfrak{B}^* \vdash_{\texttt{IL}} A$ holds via $\mathscr{D}_2^*$, and $C, \Gamma, \Re^*, \mathfrak{B}^* \vdash_{\texttt{IL}} A$ holds via $\mathscr{D}_3^*$. But then, $\Gamma, \Re^*, \mathfrak{B}^* \vdash_{\texttt{IL}} A$ holds too, via
            \begin{prooftree}
                \AxiomC{$\Gamma, \Re^*, \mathfrak{B}^*$}
                \noLine
                \UnaryInfC{$\mathscr{D}^*_1$}
                \noLine
                \UnaryInfC{$B \vee C$}
                \AxiomC{$[B], \Gamma, \Re^*, \mathfrak{B}^*$}
                \noLine
                \UnaryInfC{$\mathscr{D}_2^*$}
                \noLine
                \UnaryInfC{$A$}
                \AxiomC{$[C], \Gamma, \Re^*, \mathfrak{B}^*$}
                \noLine
                \UnaryInfC{$\mathscr{D}_3^*$}
                \noLine
                \UnaryInfC{$A$}
                \TrinaryInfC{$A$}
            \end{prooftree}
            \item application of ($\wedge_I$), ($\wedge_E$), ($\vee_I$), ($\rightarrow_I$), ($\rightarrow_E$), or ($\bot$) $\Longrightarrow$ similar to the case for ($\vee_E$);
            \item application of some rule $R$ of level $n$ in $\mathfrak{B} \Longrightarrow A \in \texttt{ATOM}_\mathscr{L}$ and $\mathscr{D}$ has the form
            \begin{prooftree}
                \AxiomC{$[\mathfrak{C}_1], \Gamma, \Re$}
                \noLine
                \UnaryInfC{$\mathscr{D}_1$}
                \noLine
                \UnaryInfC{$B_1$}
                \AxiomC{}
                \noLine
                \UnaryInfC{$\dots$}
                \AxiomC{$[\mathfrak{C}_n], \Gamma, \Re$}
                \noLine
                \UnaryInfC{$\mathscr{D}_n$}
                \noLine
                \UnaryInfC{$B_n$}
                \RightLabel{$R$}
                \TrinaryInfC{$A$}
            \end{prooftree}
            with $B_i \in \texttt{ATOM}_\mathscr{L}$ and $\mathfrak{C}_i$ set of rules of level $\leq n - 2$ discharged by $R$ ($i \leq n$). By induction hypothesis, there is $\mathscr{D}_i^*$ such that $\mathfrak{C}_i^*, \Gamma, \Re^*, \mathfrak{B}^* \vdash_{\texttt{IL}} B_i$ holds via $\mathscr{D}_i^*$ ($i \leq n$). Observe now that we have
            \begin{center}
                $R^* = \bigwedge_{i \leq n} (\bigwedge_{\rho \in \mathfrak{C}_i} \rho^* \rightarrow B_i) \rightarrow A \in \mathfrak{B}^*$.
            \end{center}
            So, $\Gamma, \Re^*, \mathfrak{B}^* \vdash_{\texttt{IL}} A$ holds too via
            \begin{prooftree}
                \AxiomC{$[\bigwedge_{\rho \in \mathfrak{C}_1} \rho^*]$}
                \UnaryInfC{$\mathfrak{C}^*_1$}
                \AxiomC{$\Gamma, \Re^*, \mathfrak{B}^*$}
                \noLine
                \BinaryInfC{$\mathscr{D}^*_1$}
                \noLine
                \UnaryInfC{$B_1$}
                \UnaryInfC{$\bigwedge_{\rho \in \mathfrak{C}_1} \rho^* \rightarrow B_1$}
                \AxiomC{}
                \noLine
                \UnaryInfC{$\dots$}
                \AxiomC{$[\bigwedge_{\rho \in \mathfrak{C}_n} \rho^*]$}
                \UnaryInfC{$\mathfrak{C}^*_n$}
                \AxiomC{$\Gamma, \Re^*, \mathfrak{B}^*$}
                \noLine
                \BinaryInfC{$\mathscr{D}^*_n$}
                \noLine
                \UnaryInfC{$B_n$}
                \UnaryInfC{$\bigwedge_{\rho \in \mathfrak{C}_n} \rho^* \rightarrow B_n$}
                \TrinaryInfC{$\bigwedge_{i \leq n} (\bigwedge_{\rho \in \mathfrak{C}_i} \rho^* \rightarrow B_i)$}
                \AxiomC{$\bigwedge_{i \leq n} (\bigwedge_{\rho \in \mathfrak{C}_i} \rho^* \rightarrow B_i) \rightarrow A$}
                \BinaryInfC{$A$}
            \end{prooftree}
            \item application of some rule $R$ of level $n$ in $\Re \Longrightarrow$ similar to the previous case, except that now we have 
            \begin{center}
                $R^* = \bigwedge_{i \leq n} (\bigwedge_{\rho \in \mathfrak{C}_i} \rho^* \rightarrow B_i) \rightarrow A \in \Re^*$.
            \end{center}
        \end{itemize}
    \end{itemize}
    ($\Longleftarrow$) We proceed by induction on the length $\lambda(\mathscr{D})$ of the derivation $\mathscr{D}$ in $\texttt{IL}$:
    \begin{itemize}
        \item $\lambda(\mathscr{D}) = 0 \Longrightarrow$ there are three cases. First, $A \in \Gamma$, hence $\Gamma \vdash_{\texttt{IL} \cup \mathfrak{B}} A$, hence $\Gamma, \Re \vdash_{\texttt{IL} \cup \mathfrak{B}} A$. The other two cases are $A \in \Re^*$ or $A \in \mathfrak{B}^*$. We proceed for both by induction on the level of the rule $R_A$ such that $R_A^* = A$. Suppose first that $A \in \Re^*$. We show $R_A \vdash_{\texttt{IL} \cup \mathfrak{B}} A$. If $R_A \in \Re$ has level $0$, the result holds trivially. Suppose the theorem proved for all rules of level $\leq n$, and suppose that $R_A$ has level $n + 1$. Hence, $R_A$ has the form
        \begin{prooftree}
            \AxiomC{$[\mathfrak{C}_1]$}
            \noLine
            \UnaryInfC{$B_1$}
            \AxiomC{}
            \noLine
            \UnaryInfC{$\dots$}
            \AxiomC{$[\mathfrak{C}_m]$}
            \noLine
            \UnaryInfC{$B_m$}
            \RightLabel{$R_A$}
            \TrinaryInfC{$B$}
        \end{prooftree}
        where $\mathfrak{C}_i$ is a set of rules of level $\leq n - 1$ discharged by $R_A$ ($i \leq m$). So,
        \begin{center}
            $A = R_A^* = \bigwedge_{i \leq m}(\bigwedge_{\rho \in \mathfrak{C}_i} \rho^* \rightarrow B_i) \rightarrow B$.
        \end{center}
        By induction hypothesis, for every $i \leq m$ and every $\rho \in \mathfrak{C}_i$, $\rho \vdash_{\texttt{IL} \cup \mathfrak{B}} \rho^*$. This means that we have a $\mathscr{D}_i$ in $\texttt{IL} \cup \mathfrak{B}$ proving $\mathfrak{C}_i \vdash_{\texttt{IL} \cup \mathfrak{B}} \bigwedge_{\rho \in \mathfrak{C}_i} \rho^*$. But then we can build
        \begin{prooftree}
            \AxiomC{$[\mathfrak{C}_1]$}
            \noLine
            \UnaryInfC{$\mathscr{D}_1$}
            \noLine
            \UnaryInfC{$\bigwedge_{\rho \in \mathfrak{C}_1} \rho^*$}
            \AxiomC{$[\bigwedge_{i \leq m}(\bigwedge_{\rho \in \mathfrak{C}_i} \rho^* \rightarrow B_i)]$}
            \UnaryInfC{$\bigwedge_{\rho \in \mathfrak{C}_1} \rho^* \rightarrow B_1$}
            \BinaryInfC{$B_1$}
            \AxiomC{}
            \noLine
            \UnaryInfC{$\dots$}
            \AxiomC{$[\mathfrak{C}_m]$}
            \noLine
            \UnaryInfC{$\mathscr{D}_m$}
            \noLine
            \UnaryInfC{$\bigwedge_{\rho \in \mathfrak{C}_m} \rho^*$}
            \AxiomC{$[\bigwedge_{i \leq m}(\bigwedge_{\rho \in \mathfrak{C}_i} \rho^* \rightarrow B_i)]$}
            \UnaryInfC{$\bigwedge_{\rho \in \mathfrak{C}_m} \rho^* \rightarrow B_m$}
            \BinaryInfC{$B_m$}
            \RightLabel{$R_A$}
            \TrinaryInfC{$B$}
            \UnaryInfC{$\bigwedge_{i \leq m}(\bigwedge_{\rho \in \mathfrak{C}_i} \rho^* \rightarrow B_i) \rightarrow B$}
        \end{prooftree}
        whence $R_A \vdash_{\texttt{IL} \cup \mathfrak{B}} A$. Since $R_A \in \Re$, hence, $\Gamma, \Re \vdash_{\texttt{IL} \cup \mathfrak{B}} A$. As regards the case when $A \in \mathfrak{B}^*$, the reasoning is the same as in the previous case, except that rules have not to be assumed, but drawn from $\mathfrak{B}$, so we obtain $\vdash_{\texttt{IL} \cup \mathfrak{B}} A$, hence $\Gamma, \Re \vdash_{\texttt{IL} \cup \mathfrak{B}} A$.
    \end{itemize}
    When $\lambda(\mathscr{D}) = k + 1$, $\mathscr{D}$ will end by applying some rule from $\texttt{IL}$, so induction can proceed smoothly.
\end{proof}

\begin{corollary}
    $\Gamma \vdash_{\emph{\texttt{IL}} \cup \mathfrak{B}} A \Longleftrightarrow \exists \Delta \subseteq \mathfrak{B}^*$ finite such that $\Gamma, \Delta \vdash_{\emph{\texttt{IL}}} A$
\end{corollary}

\begin{theorem}
   $\emph{\texttt{IL}}$ is base-complete over $\Vdash_n \ \Longleftrightarrow \ \Vdash_n$ enjoys the compact export principle and $\emph{\texttt{IL}}$ is complete over $\Vdash_n$.
\end{theorem}

\begin{proof}
    ($\Longrightarrow$) That base-completeness generally implies completeness has been already established in Proposition 12. Let us show that base-completeness of $\texttt{IL}$ over $\Vdash_n$ implies that the latter enjoys the compact export principle. Suppose $\Gamma \Vdash_{\mathfrak{B}, \ n} A$. By base-completeness, $\Gamma \vdash_{\texttt{IL} \cup \mathfrak{B}} A$ and, by Corollary 12, $\exists \Delta \subseteq \mathfrak{B}^*$ finite such that $\Gamma, \Delta \vdash_{\texttt{IL}} A$. By soundness of $\texttt{IL}$, then, $\exists \Delta \subseteq \mathfrak{B}^*$ finite such that $\Gamma, \Delta \Vdash_n A$. Next, suppose $\exists \Delta \subseteq \mathfrak{B}^*$ finite such that $\Gamma, \Delta \Vdash_n A$. Since we assumed base-completeness, which implies completeness, we have $\exists \Delta \subseteq \mathfrak{B}^*$ finite such that $\Gamma, \Delta \vdash_{\texttt{IL}} A$ which, by Corollary 12, implies $\Gamma \vdash_{\texttt{IL} \cup \mathfrak{B}} A$. By Proposition 14, i.e., base-soundness of $\texttt{IL}$, we have hence $\Gamma \Vdash_{\mathfrak{B}, \ n} A$. ($\Longleftarrow$) Suppose $\Gamma \Vdash_{\mathfrak{B}, \ n} A$. Since we are assuming that $\Vdash_n$ enjoys the compact export principle, we have that $\exists \Delta \subseteq \mathfrak{B}^*$ finite such that $\Gamma, \Delta \Vdash_n A$ and, since we are assuming completeness of $\texttt{IL}$ over $\Vdash_n$, we have that $\exists \Delta \subseteq \mathfrak{B}^*$ finite such that $\Gamma, \Delta \vdash_{\texttt{IL}} A$. Then, by Corollary 12, we have $\Gamma \vdash_{\texttt{IL} \cup \mathfrak{B}} A$.
\end{proof}

\noindent Theorem 18 plus Corollary 11 imply that base-completeness of $\texttt{IL}$ is an inconsistent notion. From which it also follows that $\texttt{IL}$ is not base-complete over \emph{any} of the proof-theoretic semantic approaches we have discussed so far, for \emph{any} $n$. We have already proved this for $\models_n$ and $\models^\alpha_n$ on independent grounds in Proposition 13, and now we can establish the same for Sandqvist's approach as a corollary of Theorem 18.

\begin{corollary}
    $\emph{\texttt{IL}}$ is not base-complete over $\models^s_n$. 
\end{corollary}

\noindent Observe that, since our proofs above were independent of the fact that atomic bases are of an intuitionistic kind, they also hold in the case when they are not, i.e., they also apply to Sandqvist's framework properly understood---that is, where $\bot$ is understood as a nullary operator defined by a special clause, and by taking atomic bases which are not necessarily required to contain the atomic explosion rules.

Let me just remark that base-incompleteness of $\texttt{IL}$ over $\models^s_2$ might have been proved more directly through the observation that $GDP^2$ fails over Sandqvist's consequence with atomic bases of level $\geq 2$, whereas Corollary 12 implies the following.
    
    \begin{proposition}
        With $\Gamma$ disjunction free, $\Gamma \vdash_{\emph{\texttt{IL}} \cup \mathfrak{B}} A \vee B \Longrightarrow \Gamma \vdash_{\emph{\texttt{IL}} \cup \mathfrak{B}} A$ or $\Gamma \vdash_{\emph{\texttt{IL}} \cup \mathfrak{B}} B$.
    \end{proposition}

\noindent Another quicker proof of the same obtains by combining Proposition 16 and Corollary 12. Be that as it may, it seems to me to have a broader interest the fact that, for every $\Sigma$ enjoying the conditions of Theorem 19, base-completeness of $\Sigma$ over $\Vdash_n$ is equivalent to compact export principle of $\Sigma$ plus completeness of $\Sigma$ over $\Vdash_n$, and is hence a contradictory notion---additionally, this holds for every $n$, and not just for $n = 2$.

More in general, by Theorem 18, base-completeness of a (compact) logic $\Sigma$ over $\Vdash_n$ may fail for two reasons: either because $\Sigma$ is not complete over $\Vdash_n$---as with $\models_n$ and $\models^\alpha_n$---or because $\Vdash_n$ is not compact at the level of consequence over atomic bases, while base-completeness forces such a ‘‘local" kind of compactness---as with $\models^s_2$. This is what the following result states.

\begin{theorem}
   $\Sigma$ is base-sound over $\Vdash_n$, enjoys the compact export principle and does not derive Harrop's rule $\Longrightarrow \Sigma$ is not base-complete over $\Vdash_n$.
\end{theorem}

\section{Concluding remarks}

The crucial incompleteness results obtained by de Campos Sanz, Piecha and Schroeder-Heister \cite{piecha, piechaschroeder-heisterdecampossanz, piechaschroeder-heister2019} are referred to a broad framework where, either atomic bases are ordered by inclusion but without requiring them to have an upper-bounded level, or else no constraint is put on the structure of the underlying set of atomic bases. The same applies to more recent findings \cite{schroederheisterrolf, stafford1, stafford2}. The results presented in this paper are on the contrary referred to sets of atomic bases where the extension-relation is given in terms of inclusion, and where the rules-level must not be greater than a fixed bound. The last constraint was forced by the fact that, among the aims of the paper, there was that of applying the results of Section 4 to (in)completeness issues, which in turn required in some cases---such as Sandqvist's---to focus on atomic rules of level 2. The upper bound limitation can be however dropped out, while retaining most of the results presented in this paper.\footnote{As done in Section 4, I shall indicate these unlimited notions by $\models$, $\models^\alpha$ and $\models^s$ respectively, and I shall use the notation $\Vdash$ to refer to any of them. So, we have what follows---I will just limit myself to the main results, but the reader should understand these as given in a context where all the other notions come with the limited consequence relations replaced by their unlimited counter-parts.

\begin{theorem}
    $\Gamma \models_\mathfrak{B} A \Longleftrightarrow \Gamma \models^\alpha_\mathfrak{B} A$.
\end{theorem}

\begin{theorem}
    $\Gamma \models A \Longleftrightarrow \Gamma \models^\alpha A$
\end{theorem}

\begin{theorem}
    $\exists \Gamma \ \exists A \ (\Gamma \models^\alpha A$ and $\Gamma \not\vdash_{\emph{\texttt{IL}}} A)$.
\end{theorem}

\begin{theorem}
    $\forall \Gamma \ \forall A \ \forall \mathfrak{C} \supseteq \mathfrak{B} \ (\Gamma \models^s_\mathfrak{C} A \Longrightarrow \Gamma \models^\alpha_\mathfrak{C} A) \Longrightarrow \ \forall \Gamma \ \forall A \ \forall \mathfrak{C} \supseteq \mathfrak{B} \ (\Gamma \models^\alpha_\mathfrak{C} A \Longrightarrow \Gamma \models^s_\mathfrak{C} A)$
\end{theorem}

\begin{corollary}
    $\forall \Gamma \ \forall A \ \forall \mathfrak{B} \in \mathbb{B} \ (\Gamma \models^s_\mathfrak{B} A \Longrightarrow \Gamma \models^\alpha_\mathfrak{B} A) \Longrightarrow \ \forall \Gamma \ \forall A \ \forall \mathfrak{B} \in \mathbb{B} \ (\Gamma \models^\alpha_\mathfrak{B} A \Longrightarrow \Gamma \models^s_\mathfrak{B} A)$
\end{corollary}

\begin{theorem}
    Both the antecedent and the consequent of Corollary 15 fail.
\end{theorem}

\begin{theorem}
    $\exists \Sigma \ (\Sigma$ base-complete over $\models^s$ and base-sound over $\models^\alpha) \ \Longrightarrow \ \forall \Gamma \ \forall A \ (\Gamma \models^s A \Longleftrightarrow \Gamma \models^\alpha A)$.
\end{theorem}

\begin{corollary}
    $\Sigma$ base-complete over $\models^s$ and base-sound over $\models^\alpha \ \Longrightarrow \ \Sigma$ complete over $\models^\alpha$.
\end{corollary}

\begin{proposition}
    $\forall \mathfrak{B} \in \mathbb{B} \ (\Gamma, \Re \vdash_{\emph{\texttt{IL}} \cup \mathfrak{B}} A \Longleftrightarrow \Gamma, \Re^*, \mathfrak{B}^* \vdash_\emph{\texttt{IL}} A)$.
\end{proposition}

\begin{corollary}
    $\forall \mathfrak{B} \in \mathbb{B} \ (\Gamma \vdash_{\emph{\texttt{IL}} \cup \mathfrak{B}} A \Longleftrightarrow \Gamma, \mathfrak{B}^* \vdash_\emph{\texttt{IL}} A)$.
\end{corollary}

\begin{theorem}
    $\emph{\texttt{IL}}$ is base-complete over $\Vdash \ \Longleftrightarrow \ \Vdash$ enjoys the compact export principle and $\emph{\texttt{IL}}$ is complete over $\Vdash$.
\end{theorem}

\begin{corollary}
    $\emph{\texttt{IL}}$ is not base-complete over $\Vdash$.
\end{corollary}

\begin{theorem}
    If $\Sigma$ enjoys the compact export principle and does not derive Harrop's rule, then $\Sigma$ is not base-complete over $\Vdash$.
\end{theorem}
\noindent The proofs of these results run in very much the same way as those for their limited counter-parts throughout Sections 4 and 5 in this paper.} It remains to be settled whether similar results would hold also under more liberal orders on atomic bases.

\paragraph{Acknowledgments} I am grateful to Ansten Klev, Hermógenes Oliveira, Thomas Piecha, Dag Prawitz, Peter-Schroeder-Heister, Will Stafford, and the anonymous reviewers, for precious remarks which helped me improve previous versions of this paper. This work has been supported by the grant PI 1965/1-1 for the DFG project \emph{Revolutions and paradigms in logic. The case of proof-theoretic semantics}.

\paragraph{Conflict of interests} The author declares that there is no conflict of interests.

\bibliographystyle{abbrv}
\bibliography{bibliography}

\begin{thebibliography}{10}

\bibitem{ayhangoodreductions}
S.~Ayhan.
\newblock What are acceptable reductions? {Perspectives} from proof-theoretic
  semantics and type theory.
\newblock {\em Australasian journal of logic}, 20(3):412--428, 2023.
\newblock \url{https://doi.org/10.26686/ajl.v20i3.7692}.

\bibitem{barrosopereiranascimento}
V.~{Barroso Nascimento}, L.~Pereira, and E.~Pimentel.
\newblock An ecumenical view of proof-theoretic semantics.
\newblock {\em Synthese}, 2025.
\newblock \url{https://doi.org/10.1007/s11229-025-05269-z}.

\bibitem{piechadecampossanz}
W.~{de Campos Sanz} and T.~Piecha.
\newblock A critical remark on the {BHK} interpretation of implication.
\newblock {\em Philosophia Scientiae}, 18(3):13--22, 2014.
\newblock \url{https://doi.org/10.4000/philosophiascientiae.965}.

\bibitem{gentzen}
G.~Gentzen.
\newblock Untersuchungen über das logische {Schließen} i, ii.
\newblock {\em Mathematische Zeitschrift}, 39:176–210, 405–431, 1935.
\newblock \url{https://doi.org/10.1007/BF01201353}.

\bibitem{gheorghiupymnegationasfailure}
A.~Gheorghiu and D.~Pym.
\newblock Definite formulae, negation-as-failure, and the base-extension
  semantics of intuitionistic propositional logic.
\newblock {\em Bulletin of the section of logic}, 52(2):239--266, 2023.
\newblock \url{https://doi.org/10.18778/0138-0680.2023.16}.

\bibitem{GheorghiuPym}
A.~Gheorghiu and D.~Pym.
\newblock From proof-theoretic validity to base-extension semantics for
  intuitionistic propositional logic.
\newblock {\em Studia Logica}, 2025.
\newblock \url{https://doi.org/10.1007/s11225-024-10163-9}.

\bibitem{harrop}
R.~Harrop.
\newblock Concerning formulas of the types ${A} \rightarrow {B} \vee {C}$, ${A}
  \rightarrow (\exists x){B}(x)$ in intuitionistic formal systems.
\newblock {\em Journal of symbolic logic}, 25:27--32, 1960.
\newblock \url{doi:10.2307/2964334}.

\bibitem{howard}
W.~Howard.
\newblock The formulae-as-types notion of construction.
\newblock In J.~Seldin and R.~Hindley, editors, {\em To {H.B. Curry}: Essays on
  combinatory logic, lambda calculus and formalism}, pages 479--490. Academy
  Press, 1980.

\bibitem{martin-loef}
P.~{Martin-L\"{o}f}.
\newblock {\em Intuitionistic type theory}.
\newblock Bibliopolis, 1984.

\bibitem{hermogenespragmatist}
H.~Oliveira.
\newblock On {Dummett's} pragmatist justification procedure.
\newblock {\em Erkenntnis}, 86:429--455, 2021.
\newblock \url{https://doi.org/10.1007/s10670-019-00112-7}.

\bibitem{piccolominibook}
A.~{Piccolomini d'Aragona}.
\newblock {\em Prawitz's epistemic grounding. An investigation into the power
  of deduction}.
\newblock Synthese Library, Springer, 2023.

\bibitem{ptssquare}
A.~{Piccolomini d'Aragona}.
\newblock The proof-theoretic square.
\newblock {\em Synthese}, 201(219), 2023.
\newblock \url{https://doi.org/10.1007/s11229-023-04203-5}.

\bibitem{piccolominiattitudes}
A.~{Piccolomini d'Aragona}.
\newblock Four constructivist attitudes in {Prawitzian} semantics.
\newblock In I.~Sedlar, editor, {\em The Logica Yearbook 2022}. College
  Publications, 2024.

\bibitem{piccolomininote}
A.~{Piccolomini d'Aragona}.
\newblock A note on schematic validity and completeness in {Prawitz}'s
  semantics.
\newblock In F.~Bianchini, V.~Fano, and P.~Graziani, editors, {\em Current
  topics in logic and the philosophy of science. Papers from SILFS 2022
  postgraduate conference}. College Publications, 2024.

\bibitem{piecha}
T.~Piecha.
\newblock Completeness in proof-theoretic semantics.
\newblock In T.~Piecha and P.~Schroeder-Heister, editors, {\em Advances in
  proof-theoretic semantics}, pages 231--251. 2016.
\newblock \url{https://doi.org/10.1007/978-3-319-22686-6_15}.

\bibitem{piechaschroeder-heisterdecampossanz}
T.~Piecha, W.~{de Campos Sanz}, and P.~Schroder-Heister.
\newblock Failure of completeness in proof-theoretic semantics.
\newblock {\em Journal of philosophical logic}, 44:321--335, 2015.
\newblock \url{https://doi.org/10.1007/s10992-014-9322-x}.

\bibitem{piechaschroeder-heister2019}
T.~Piecha and P.~{Schroeder-Heister}.
\newblock Incompleteness of intuitionistic propositional logic with respect to
  proof-theoretic semantics.
\newblock {\em Studia Logica}, 107(1):47--62, 2019.
\newblock \url{https://doi.org/10.1007/s11225-018-9823-7}.

\bibitem{prawitz1965}
D.~Prawitz.
\newblock {\em Natural deduction. A proof-theoretical study}.
\newblock Almqvist \& Wiskell, 1965.

\bibitem{prawitz1971}
D.~Prawitz.
\newblock Ideas and results in proof theory.
\newblock In J.~E. Fenstad, editor, {\em Proceedings of the second Scandinavian
  logic symposium}, pages 235--307. Elsevier, 1971.
\newblock \url{https://doi.org/10.1016/S0049-237X(08)70849-8}.

\bibitem{prawitz1973}
D.~Prawitz.
\newblock Towards a foundation of a general proof-theory.
\newblock In P.~Suppes, L.~Henkin, A.~Joja, and G.~C. Moisil, editors, {\em
  Proceedings of the Fourth International Congress for Logic, Methodology and
  Philosophy of Science, Bucharest, 1971}, pages 225--250. Elsevier, 1973.
\newblock \url{https://doi.org/10.1016/S0049-237X(09)70361-1}.

\bibitem{prawitz2015}
D.~Prawitz.
\newblock Explaining deductive inference.
\newblock In H.~Wansing, editor, {\em {Dag Prawitz} on proofs and meaning},
  pages 65--100. Springer, 2015.

\bibitem{sandqvist}
T.~Sandqvist.
\newblock Base-extension semantics for intuitionistic sentential logic.
\newblock {\em Logic journal of the IGPL}, 23(5):719--731, 2015.
\newblock \url{https://doi.org/10.1093/jigpal/jzv021}.

\bibitem{schroederheisternaturalext}
P.~Schroeder-Heister.
\newblock A natural extension of natural deduction.
\newblock {\em The journal of symbolic logic}, 49(4):1284--1300, 1984.
\newblock \url{https://doi.org/10.2307/2274279}.

\bibitem{schroeder-heister2006}
P.~Schroeder-Heister.
\newblock Validity concepts in proof-theoretic semantics.
\newblock {\em Synthese}, 148:525--571, 2006.
\newblock \url{https://doi.org/10.1007/s11229-004-6296-1}.

\bibitem{schroeder-heisterSE}
P.~Schroeder-Heister.
\newblock Proof-theoretic semantics.
\newblock In E.~N. Zalta, editor, {\em The Stanford Encyclopedia of
  Philosophy}. 2018.

\bibitem{schroederheisterrolf}
P.~Schroeder-Heister.
\newblock Prawitz's completeness conjecture: a reassessment.
\newblock {\em Theoria}, 90(5):515--527, 2024.
\newblock \url{https://doi.org/10.1111/theo.12541}.

\bibitem{stafford1}
W.~Stafford.
\newblock Proof-theoretic semantics and inquisitive logic.
\newblock {\em Journal of philosophical logic}, 50:1199--1229, 2021.
\newblock \url{https://doi.org/10.1007/s10992-021-09596-7}.

\bibitem{stafford2}
W.~Stafford and V.~Nascimento.
\newblock Following all the rules: intuitionistic completeness for generalised
  proof-theoretic validity.
\newblock {\em Analysis}, 2023.
\newblock \url{https://doi.org/10.1093/analys/anac100}.

\bibitem{troelstravandalen}
A.~S. Troelstra and D.~van Dalen.
\newblock {\em Constructivism in mathematics. Vol. I}.
\newblock North-Holland, 1988.

\end{thebibliography}

\end{document}